\documentclass[12pt]{amsart}

\usepackage{amsthm,amsfonts,amsmath,amssymb,latexsym,epsfig}

\DeclareMathAlphabet\oldmathcal{OMS}        {cmsy}{b}{n}
\SetMathAlphabet    \oldmathcal{normal}{OMS}{cmsy}{m}{n}
\DeclareMathAlphabet\oldmathbcal{OMS}       {cmsy}{b}{n}
 
\usepackage{eucal}

\newtheorem{theorem}{Theorem}
\newtheorem{lemma}[theorem]{Lemma}
\newtheorem{proposition}[theorem]{Proposition}
\newtheorem{corollary}[theorem]{Corollary}
\newtheorem{definition}[theorem]{Definition}
\newenvironment{example}{\medskip \refstepcounter{theorem}
\noindent  {\bf Example \thetheorem}.\rm}{\,}
\newenvironment{remark}{\medskip \refstepcounter{theorem}
\noindent  {\bf Remark \thetheorem}.\rm}{\,}

\renewcommand{\thetheorem}{\thesection.\arabic{theorem}}
\def\d{\partial}                                   
\def\db{\overline{\partial}}
\def\ddb{\partial  \overline{\partial}}  
\def\<{\langle}

\def\>{\rangle}                                    
\def\a{\alpha}                                     
\def\b{\beta}
\def\g{\gamma}                                     
\def\v{\varphi}                                    
                                     
\def\o{\omega}

\def\hok{\mbox{}\begin{picture}(10,10)\put(1,0){\line(1,0){7}}
  \put(8,0){\line(0,1){7}}\end{picture}\mbox{}}
\def\BOne{{\mathchoice {\rm 1\mskip-4mu l} {\rm 1\mskip-4mu l}
                          {\rm 1\mskip-4.5mu l} {\rm 1\mskip-5mu l}}}
\def\Ric{{\rm Ric}}

\def\fract#1#2{\raise4pt\hbox{$ #1 \atop #2 $}}

\def\bbr{{\mathbb R}}

\def\bbu{{\mathbb U}}

\def\grz{\zeta}
\def\gro{\omega}

\def\gra{\alpha}
\def\grd{\delta}

\def\cald{{\mathcal D}}

\def\calf{{\mathcal F}}

\def\cals{{\oldmathcal S}}

\def\calx{{\mathcal X}}

\def\la#1{\hbox to #1pc{\leftarrowfill}}
\def\ra#1{\hbox to #1pc{\rightarrowfill}}

\def\ga{{\mathfrak a}}

\def\gc{{\mathfrak c}}

\def\gf{{\mathfrak f}}

\def\gh{{\mathfrak h}}
\def\gi{{\mathfrak i}}

\def\gl{{\mathfrak l}}

\def\gn{{\mathfrak n}}
\def\go{{\mathfrak o}}
\def\gp{{\mathfrak p}}

\def\gr{{\mathfrak r}}
\def\gs{{\mathfrak s}}
\def\gt{{\mathfrak t}}
\def\gu{{\mathfrak u}}

\def\gz{{\mathfrak z}}

\def\gA{{\mathfrak A}}

\def\gC{{\mathfrak C}}
\def\gD{{\mathfrak D}}

\def\gF{{\mathfrak F}}

\def\gH{{\mathfrak H}}

\def\gR{{\mathfrak R}}
\def\gS{{\mathfrak S}}
\def\gT{{\mathfrak T}}

\newcommand{\CX}{\mbox{${\calx}\hspace{-.8em}-\,$}}

\begin{document}

\title[Canonical Sasakian metrics]{Canonical Sasakian metrics}

\author[C.P. Boyer]{Charles P. Boyer}
\author[K. Galicki]{Krzysztof Galicki}
\author[S.R. Simanca]{Santiago R. Simanca}

\thanks{During the preparation of this work, the first two authors
were partially supported by NSF grant DMS-0504367.}

\address{
Department of Mathematics and Statistics, 
University of New Mexico, Albuquerque, N.M. 87131}
\email{cboyer@math.unm.edu, galicki@math.unm.edu, santiago@math.unm.edu}

\begin{abstract}
Let $M$ be a closed manifold of Sasaki type. A polarization of $M$ is defined
by a Reeb vector field, and for any such polarization, we consider the set of 
all 
Sasakian metrics compatible with it. On this space we study the functional 
given by the square of the $L^2$-norm of the scalar curvature. We prove 
that its critical points, or canonical representatives of the polarization, 
are Sasakian metrics that are transversally extremal. We define a Sasaki-Futaki
invariant of the polarization, and show that it obstructs the existence of 
constant scalar curvature representatives. For a fixed CR structure of 
Sasaki type, we define the Sasaki cone of structures compatible with this 
underlying CR structure, and prove that the set of polarizations in it that 
admit a canonical representative is open. We use our results to describe
fully the case of the sphere with its standard CR structure, showing that each
element of its Sasaki cone can be represented by a canonical metric; we
compute their Sasaki-Futaki invariant, and use it to describe the canonical
metrics that have constant scalar curvature, and to prove that just the 
standard polarization can be represented by a Sasaki-Einstein metric.
\end{abstract}

\maketitle

\section{Introduction}
With the knowledge that the set of K\"ahler metrics representing a
given K\"ahler class is an affine space modeled after the smooth functions,
Calabi used \cite{cal0,cal1} a natural Riemannian functional on this space
with the hope of using it to find canonical representatives of the given class.
In effect, his functional, or Calabi energy, is simply the squared $L^2$-norm 
of the 
scalar curvature, and the critical point minimizing it would fix the affine 
parameter alluded to above, yielding the desired representative of the 
class. Calabi named these critical points extremal K\"ahler metrics.
It was then determined that if the Futaki character \cite{fu} of the class 
vanishes, a plausible extremal K\"ahler representative must be a metric of 
constant scalar curvature, and if under that condition we look at the case
where the K\"ahler class in question is a multiple of the first Chern class, 
the extremal representative must then be K\"ahler-Einstein.

One of the most important problems in K\"ahler geometry
today involves the subtle questions regarding the existence of
extremal K\"ahler metrics representing a given class.  Over the years,
starting with the formulation of the famous Calabi
Conjecture and its proof by Yau in 1978, various tools have been used
or developed to attack this problem. The continuity method, Tian's
$\alpha$-invariant, the Calabi-Lichnerowicz-Matsushima obstruction,
the Futaki invariant and its generalizations, the Mabuchi K-energy,
and more recently the various notions of stability proposed and
studied by Tian, Donaldson and others. Substantial progress has
been made, but the general existence problem remains open.

Sasakian geometry sits naturally in between two K\"ahler
geometries. On the one hand, Sasakian manifolds are the bases of metric
cones which are K\"ahler. On the other hand, any Sasakian manifold
is contact, and the one dimensional foliation associated to the
characteristic Reeb vector field is transversally
K\"ahler. In many interesting situations, the orbits of the Reeb vector field 
are all closed, in which case the Sasakian structure is called quasi-regular. 
Compact quasi-regular Sasakian manifolds have the structure of
an orbifold circle bundle over a compact K\"ahler orbifold, which must be 
algebraic and which has at most cyclic quotient singularities. Since much of 
the study of compact K\"ahler manifolds and extremal metrics can be extended 
to the orbifold case, extension often done in a fairly straightforward way, it
is not surprising that we can then ``translate'' statements involving compact 
K\"ahler orbifolds to
conclude parallel statements regarding quasi-regular Sasakian
structures. This is an approach that has been spectacularly successful in
constructing new quasi-regular Sasaki-Einstein metrics on
various contact manifolds of odd dimension greater than 3
(cf. \cite{BGK05, BG06, BoGa05a, Kol05a, Kol05b}, and references therein).

For some time now, it has been believed that the only interesting
(canonical) metrics in Sasakian geometry occur precisely in this
orbibundle setting. In 1994, Cheeger and Tian conjectured that any 
compact Sasaki-Einstein manifold must be quasi-regular \cite{ChTi94}. 
Their conjecture was phrased in terms
of the properties of the Calabi-Yau cone rather than its
Sasaki-Einstein base\footnote{More precisely Cheeger and Tian
used the term {\it standard cone}, and the conjecture states that
all Calabi-Yau cones are standard \cite{ChTi94}}, and until recently,
compact Sasakian manifolds with non-closed leaves were certainly
known, but there was no evidence to suspect that we could get such
structures with Einstein metrics as well. Hence, it was reasonable
to believe that all Sasaki-Einstein metrics could be understood
well by simply studying the existence of K\"ahler-Einstein metrics on compact 
cyclic orbifolds.

As it turns out, the conjecture of Cheeger and Tian mentioned above is not 
true, and the first examples
of irregular Sasaki-Einstein manifolds, that is to say, Sasaki-Einstein
manifolds that are not quasi-regular, came first from the physics surrounding 
the famous CFT/AdS Duality Conjecture \cite{GMSW04a, GMSW04b,
MaSp05b, MaSp06, MaSpYau05, MaSpYau06, CLPP05}. It now appears
that there are irregular structures of this type on many compact
manifolds in any odd dimension greater than 3. These
Sasaki-Einstein metrics represent canonical points in the space of metrics
adapted to the underlying geometric setting. However, although their 
Calabi-Yau cones are smooth outside the tip of the cone, their space of leaves 
is not even Hausdorff, and so the whole ``orbibundle over a K\"ahler-Einstein 
base'' approach proves itself insufficient in the study of the problem.

The discovery of these new metrics make a strong case in favor of a 
variational formulation of the study of these Sasakian metrics, in a way 
analogous to the notion of the Calabi energy and extremality. With the proper 
set-up, the quasi-regularity property 
should no longer be a key factor, and all Sasaki-Einstein metrics 
should indeed appear as minima of a suitable Riemannian functional. This 
would put on equal footing the analysis of all Sasakian structures, the 
quasi-regular or the irregular ones. Thus, we should be able to study the 
existence and uniqueness of these canonical Sasakian metrics in ways 
parallel to those used in K\"ahler geometry. 

Until now, this approach for finding canonical Sasakian structures has not 
been pursued, perhaps due to the lack of evidence that the orbibundle 
approach would be insufficient. We propose here to look at the squared 
$L^2$-norm of the scalar curvature functional, defined over a suitable space 
of Sasakian metrics that are determined by fixing the Reeb vector field,
which we think of as polarizing the Sasakian manifold. Its critical points
are, by definition, canonical Sasakian metrics representing the said 
polarization. 

Recently, Martelli, Sparks and Yau presented a similar point of view 
\cite{MaSpYau06}, opting to look at the Sasaki-Einstein metrics as minima of 
the Hilbert action instead. Our point of view has important advantages, several
of which are elucidated in the present article. In particular, we see that
for certain manifolds of Sasaki type, the {\it optimal} Sasaki metric on it 
cannot possibly have constant scalar curvature, showing the need to enlarge
the plausible set of metrics to be considered, if at least, from a 
mathematical point of view.

In general, the minimization of the $L^2$-norm of the scalar curvature over 
metrics of fixed volume is intimately related to the search for Einstein 
metrics. We use this functional here over a smaller space of metrics, 
thus laying the foundation for the study of canonical Sasakian metrics 
in a way that parallels what is done in K\"ahler geometry for the extremal 
metric problem. This point of view eliminates the need to make a distinction 
between the quasi-regular and the irregular case, discussing them both 
on an equal footing. For we introduce the notion of a polarized 
Sasakian manifold, polarized by a Reeb vector field, and analyze the 
variational problem for the $L^2$-norm of the scalar curvature 
over the space of Sasakian metrics representing the said polarization. Given
a CR structure of Sasaki type, we define the cone
of Sasakian polarizations compatible with this underlying CR structure,
and discuss the variational problem for this functional as we vary the 
polarization on this cone also. The quasi-regularity or not of the
resulting critical Sasakian structures is just a property of the 
characteristic foliation defined by the Reeb vector field. This foliation 
must clearly sit well with the Sasakian metrics under consideration in our
approach, but it stands in its own right. A canonical Sasakian metric,
a critical point of the said functional,  
interacts with the underlying characteristic foliation, but neither one of them
determines the other.

We organize the paper as follows. In \S2 we recall and review the
necessary definitions of Sasakian manifolds and associated 
structures. In \S3 we define the notion of a polarized Sasakian manifold, and
describe the space of Sasakian metrics that represent a given polarization, a
space consisting of metrics with the same transversal holomorphic structure.
We then analyze the variational problem for the $L^2$-norm of the scalar 
curvature with it as its domain of definition, and show that the resulting 
critical points are Sasakian metrics for which the basic vector field
$\d_g^{\#}s_g^T=\d_g^{\#}s_g$ is transversally holomorphic, that is to say,
metrics that are transversally extremal. 
In \S4 we study various transformation groups of Sasakian structures and 
their Lie algebras, proving the Sasakian version of the Lichnerowicz-Matsushima
theorem. In \S5 we define and study the Sasaki-Futaki
invariant, and prove that a canonical Sasakian metric is of constant scalar 
curvature if, and only if, this invariant vanishes for the polarization under
consideration.
In \S6 we define and study the Sasaki cone, and end up in \S7 by proving that 
the polarizations in the Sasaki cone that admit canonical representatives form
an open set, proving that the openness theorem for the extremal
cone in K\"ahler geometry \cite{cs} holds in the Sasakian context also.
We illustrate the power of this result by providing a detailed analysis of
the Sasaki cone for the standard CR-structure on the unit sphere
${\mathbb S}^{2n+1}$, and use it to show that all of its elements admit 
canonical representatives. We describe explicitly those that are of constant 
scalar curvature, and show that the standard metric is the only one of 
these that is Sasaki-Einstein.

\section{Sasakian manifolds}
\setcounter{theorem}{0}
We recall that an almost contact structure on a differentiable manifold
$M$ is given by a triple $(\xi, \eta, \Phi)$, where $\xi$ is a vector field,
$\eta$ is a one form, and $\Phi$ is a
tensor of type $(1,1)$, subject to the relations
$$\eta(\xi) = 1 \, , \quad \Phi^2 = -\BOne + \xi \otimes \eta \, .$$
The vector field $\xi$ defines the {\it characteristic foliation} 
$\calf_{\xi}$ with one-dimensional leaves, and the kernel of $\eta$
defines the codimension one sub-bundle ${\mathcal D}$. This yields a
canonical splitting 
\begin{equation}
TM = {\mathcal D} \oplus L_{\xi}\, , \label{cs}
\end{equation}
where $L_\xi$ is the trivial line bundle generated by $\xi$.
The sub-bundle ${\mathcal D}$ inherits an almost complex structure 
$J$ by restriction of $\Phi$. Clearly, the dimension of $M$ must
be an odd integer $2n+1$. We refer to $(M,\xi,\eta, \Phi)$ as an almost
contact manifold. If we disregard the tensor $\Phi$ and
characteristic foliation, that is to say, if we just look at the sub-bundle
${\mathcal D}$ forgetting altogether its almost complex structure, we then
refer to the contact structure $(M, {\mathcal D})$, or simply
${\mathcal D}$ when $M$ is understood. Here, and further below, the reader
can no doubt observe that the historical development of the
terminology is somewhat unfortunate, and for instance, it is an 
almost contact structure the one that gives rise to a contact one, rather
than the other way around.

A Riemannian metric $g$ on $M$ is said to be compatible with the
almost contact structure $(\xi, \eta, \Phi)$ if for any pair of 
vector fields $X,Y$, we have that
$$g(\Phi(X), \Phi(Y))=g(X,Y)-\eta (X)\eta(Y) \, .$$
Any such $g$ induces an almost Hermitian metric on the sub-bundle 
${\mathcal D}$. We say that $(\xi, \eta, \Phi,g)$ is an almost contact 
metric structure.

In the presence of a compatible Riemannian metric $g$ on 
$(M,\xi, \eta, \Phi)$, the canonical decomposition (\ref{cs}) is orthogonal.
Furthermore, requiring that the orbits of the field $\xi$ be geodesics is 
equivalent to requiring that $\pounds_\xi \eta =0$, a condition that 
in view of the relation $\xi \hok \eta =1$, can be re-expressed as 
$\xi \hok d\eta =0$. 

An almost contact metric structure $(\xi,\eta,\Phi,g)$ is said to be a 
contact metric structure if for all pair of vector fields
$X$, $Y$, we have that
\begin{equation}
g(\Phi X, Y) = d\eta (X,Y)\, . \label{tkf}
\end{equation}
We then say that $(M,\xi,\eta,\Phi,g)$ is a contact metric manifold. 
Notice that in such a case, the volume element defined by $g$ is given by
\begin{equation}
d\mu_g = \frac{1}{n!}\eta \wedge (d\eta)^n  \, . \label{ve}
\end{equation}

It is convenient to reinterpret the latter structure in terms of the
cone construction. Indeed, on 
$C(M)=M\times {\mathbb R}^{+}$, we introduce the metric 
$$g_C=dr^2 + r^2 g \, .$$
The radial vector field $r\partial_r $ satisfies
$$\pounds_{r\partial_r} g_C=2g_C \, ,$$
and we may define an almost complex structure $I$ on $C(M)$ by
$$I(Y)=\Phi (Y) + \eta(Y) r\partial_r \, , \quad I(r\partial_r ) = -\xi \, .$$
The almost contact manifold $(M,\xi, \eta, \Phi)$ is said to be {\it normal} 
if the
pair $(C(M),I)$ is a complex manifold. In that case, the induced almost
complex structure $J$ on ${\mathcal D}$ is integrable.

\begin{definition}
An  contact metric structure $(\xi, \eta, \Phi, g)$ on a manifold $M$ 
is said to be a {\bf Sasakian structure} if 
$(\xi, \eta, \Phi)$ is normal. A smooth manifold provided with one such 
structure is said to be a {\bf Sasakian manifold}, or a {\bf manifold of 
Sasaki type}.
\end{definition}
\medskip

For a Sasakian structure $(\xi, \eta, \Phi, g)$, 
the integrability of the almost complex structure $I$ on
the cone $C(M)$ implies that the Reeb vector field $\xi$ leaves 
both, $\eta$ and $\Phi$, invariant \cite{BG06}. We obtain a codimension one 
integrable strictly pseudo-convex CR structure $(\cald,J)$, where 
$\cald=\ker \eta$ is the contact bundle and $J=\Phi|_\cald$, and 
the restriction of $g$ to ${\mathcal D}$ defines a positive definite 
symmetric form on $({\mathcal D},J)$ that we shall refer to as the 
transverse K\"ahler metric $g^T$. 

By (\ref{tkf}), the K\"ahler form of the transverse K\"ahler metric 
is given by the form $d\eta$. Therefore, the Sasakian metric $g$ is 
determined fully in terms 
of $(\xi,\eta, \Phi)$ by the expression
\begin{equation}
g=d\eta \circ (\BOne \otimes \Phi)+ \eta \otimes \eta \, ,\label{me1}
\end{equation}
where the fact that $d\eta$ is non-degenerate over ${\mathcal D}$ is
already built in. Since $\xi$ leaves invariant $\eta$ and $\Phi$, it is 
a Killing field, its orbits are geodesics, and the 
decomposition (\ref{cs}) is orthogonal. Despite its dependence on the 
other elements of the structure, we insist on explicitly referring to $g$ as 
part of the Sasakian structure $(\xi, \eta, \Phi, g)$. 

The discussion may be turned around to produce an alternative definition of
the notion of Sasakian structure. For the contact metric structure 
$(\xi,\eta,\Phi,g)$
defines a Sasakian structure if $({\mathcal D},J)$ is a complex sub-bundle
of $TM$, and $\xi$ generates a group of isometries. This alternative
approach appears often in the literature.

If we look at the Sasakian structure $(\xi, \eta, \Phi, g)$ from the 
point of view of CR geometry, its underlying strictly pseudo-convex CR 
structure $(\cald,J)$, with associated contact bundle ${\mathcal D}$, has
Levi form $d\eta$. (In the sequel, when 
referring to any CR structure, we shall always mean one that is integrable and
of codimension one.)

\begin{definition}\label{CRSas}
Let $(\cald,J)$ be a strictly pseudo-convex {\rm CR} structure on $M$. We say 
that $(\cald,J)$ is of {\bf Sasaki type} if there exists a Sasakian structure 
$\cals=(\xi,\eta,\Phi,g)$ such that $\cald=\ker\, \eta$ and $\Phi|_\cald=J$.
\end{definition}

If $(\cald,J)$ is a CR structure of Sasaki type, the Sasakian structures 
$\cals=(\xi,\eta,\Phi,g)$ with  $\cald=\ker\, \eta$ and $\Phi|_\cald=J$ will be
said to be Sasakian structures {\it with underlying} CR {\it structure} 
$(\cald,J)$.

The following well known result will be needed later. Let us observe 
that since the fibers of the Riemannian foliation defined by a Sasaki 
structure $(\xi,\eta,\Phi,g)$ are geodesics, their second
fundamental forms are trivial.

\begin{proposition}\label{pr1}
Let $(M,\xi,\eta,\Phi,g)$ be a Sasakian manifold. Then we have that
\begin{enumerate}
\item[a)] $Ric_g(X,\xi)=2n\eta(X)$ for any vector field $X$.
\item[b)] $Ric_g(X,Y)=Ric_T (X,Y)-2g(X,Y)$ for any pair of sections
$X,Y$ of ${\mathcal D}$.
\item[c)] $s_g=s_T-|A|^2=s_T-2n$, where $A$ is the O'Neill tensor of
the ``corresponding'' Riemannian submersion.  
\end{enumerate}
In these statements, the subscript $T$ denotes the corresponding intrinsic 
geometric quantities of the transversal metric $g^T$.
\end{proposition}

{\it Proof}. The Riemannian submersion of the Sasakian structure has totally 
geodesic fibers. 
The vector field $\xi$ spans the only vertical direction, and we have that
$A_XY=-(d\eta(X,Y))\xi$, and that $A_X \xi =-\Phi (X)$ for arbitrary 
horizontal vector fields $X,Y$. The first two results follow easily from 
O'Neill's formulae \cite{on}. The 
computation of the $L^2$-norm of $A$ is a simple consequence of the
$J$-invariance of the induced metric on ${\mathcal D}$. \qed

\section{Canonical representatives of polarized Sasakian manifolds}
\setcounter{theorem}{0} 
Let us consider a Sasakian structure $(\xi,\eta,\Phi,g)$ on $M$ that shall
remain fixed throughout this section. There are two naturally 
defined sets of deformations of this structure, those where the Reeb vector 
field remains fixed while the underlying CR structure changes by a 
diffeomorphism, and those where the underlying CR structure stays put while 
the Reeb vector field varies. We study the first type of deformations in this 
section. They turn out to be the Sasakian analogues of the
set of K\"ahler metrics representing a given polarization on a manifold 
of K\"ahler type. They have different but isomorphic underlying CR structure,
and they all share the same transverse holomorphic structure. The other
set of deformations that fix the underlying CR structure will lead to the 
Sasakian analogue of the K\"ahler cone of a manifold of K\"ahler type, and
they will be analyzed in \S\ref{sc} below.

We begin by recalling that a function $\v\in C^{\infty}(M)$ is said to be 
basic if it is annihilated by the vector field $\xi$, that is to say, if 
$\xi(\v)=0$. We denote by $C_B^{\infty}(M)$ the space of all real valued
basic functions on $M$. We observe that the notion of basic can be extended 
to covariant 
tensors of any order in the obvious manner. In particular, when looking at the 
transversal K\"ahler metric of $(\xi,\eta,\Phi,g)$, its K\"ahler form is 
basic, and so must be all of its curvature tensors as well. This observation 
will play a crucial r\^ole in the sequel. 

We consider the set \cite{BG06}
\begin{equation}
\label{Sasakianspace}
{\mathcal S}(\xi)=\{\text{Sasakian structure $(\tilde{\xi}, \tilde{\eta},
\tilde{\Phi}, \tilde{g})$} \, | \; \tilde{\xi}=\xi\}\, ,
\end{equation}
and provide it with the $C^\infty$ compact-open topology as sections of 
vector bundles.
For any element $(\tilde{\xi}, \tilde{\eta},
\tilde{\Phi}, \tilde{g})$ in this set, the 1-form $\grz=\tilde{\eta} - \eta$ 
is basic, and so $[d\tilde{\eta}]_B=[d\eta]_B$. Here, $[\, \cdot \, ]_B$ stands
for a cohomology class in the basic cohomology ring, ring that is defined by 
the restriction $d_B$ of the exterior derivative $d$ to the subcomplex of 
basic forms in the de Rham complex of $M$.
Thus, all of the Sasakian structures in ${\mathcal S}(\xi)$ correspond to the 
same  basic cohomology class. We call 
${\mathcal S}(\xi)$ the {\it space of Sasakian structures compatible with 
$\xi$}, and say that the Reeb vector field $\xi$ {\it polarizes} the 
Sasakian manifold $M$.

Given the Reeb vector field $\xi$, we have its characteristic foliation 
$\calf_\xi,$ so we let 
$\nu(\calf_\xi)$ be the vector 
bundle whose fiber at a point $p\in M$ is the quotient space $T_pM/L_{\xi}$, 
and let $\pi_\nu: TM \rightarrow \nu(\calf_\xi)$ be 
the natural projection. The background structure $\cals=(\xi,\eta,\Phi,g)$
induces a complex structure $\bar{J}$ on $\nu(\calf_\xi)$. This is
defined by $\bar{J}\bar{X}:=\overline{\Phi(X)}$, where $X$ is any 
vector field
in $M$ such that $\pi (X)=\bar{X}$. Furthermore, the underlying CR 
structure $({\mathcal D},J)$ of $\cals$ is isomorphic to
$(\nu(\calf_\xi),\bar{J})$ as a complex vector bundle. For this reason,
we refer to $(\nu(\calf_\xi),\bar{J})$ as the complex normal bundle
of the Reeb vector field $\xi$, although its identification with
$({\mathcal D},J)$ is not canonical. We shall 
say that $(M,\xi,\nu(\calf_\xi),\bar{J})$, or simply $(M,\xi, \bar{J})$, is a 
polarized Sasakian manifold.

We define ${\mathcal S}(\xi,\bar{J})$ to be the subset of all 
structures $(\tilde{\xi}, \tilde{\eta},
\tilde{\Phi}, \tilde{g})$ in ${\mathcal S}(\xi)$
such that the diagram
\begin{equation}\label{Sasakianspace2eqn}
\begin{array}{ccccc}
TM & \stackrel{\tilde{\Phi}}{\rightarrow} &  TM & \\
  \downarrow \! \! \mbox{{\small $\pi_\nu$}}  & & \downarrow \! \! 
\mbox{{\small $\pi_\nu$}} & \\
   \nu(\calf_\xi) &\stackrel{\bar{J}}{\rightarrow} &\nu(\calf_\xi),&
\end{array}
\end{equation}
commutes. This set consists of elements of ${\mathcal S}(\xi)$ 
with the same transverse holomorphic structure $\bar{J}$, or more
precisely, the same complex normal bundle $(\nu(\calf_\xi),\bar{J})$.
We have \cite{BG06}

\begin{lemma}\label{le1}
The space ${\mathcal S}(\xi, \bar{J})$ of all Sasakian structures with Reeb 
vector field $\xi$ and transverse holomorphic structure $\bar{J}$ is an affine 
space modeled on $(C_B^{\infty}(M)/{\mathbb R})\times 
(C_B^{\infty}(M)/{\mathbb R})\times H^{1}(M,{\mathbb Z})$. Indeed, if 
$(\xi, \eta, \Phi, g)$ is a given Sasakian structure in 
${\mathcal S}(\xi, \bar{J})$, any other Sasakian structure
$(\xi, \tilde{\eta},\tilde{\Phi}, \tilde{g})$ in it is determined by 
real valued basic functions $\v$ and $\psi$ and integral closed $1$-form $\a$,
such that
$$\begin{array}{rcl}
\tilde{\eta} & = & \eta+ d^c \varphi + d\psi +i(\a) \, , \\
\tilde{\Phi} & = & \Phi - (\xi \otimes (\tilde{\eta}-\eta)) \circ \Phi \, ,\\
\tilde{g} & =  & d\tilde{\eta}\circ (\BOne \otimes \tilde{\Phi} ) +
\tilde{\eta}\otimes \tilde{\eta} \, ,
\end{array}
$$ 
where $d^c=\frac{i}{2}(\db - \d)$, and $i:H^{1}(M,{\mathbb Z})\mapsto 
H^{1}(M,{\mathbb R})=H^{1}({\mathfrak F}_{\xi})$ is the homomorphism induced by
inclusion. 
In particular, $d\tilde{\eta}=d\eta + i\ddb \varphi$.
\end{lemma}

The complex structure defining the operators $\d$ and $\db$ in this Lemma is 
$\bar{J}$, as basic covariant tensors on $M$ define multilinear maps on 
$\nu({\mathcal F}_{\xi})$. We think of these as tensors on a transversal
K\"ahler manifold that does not necessarily exist. Be as it may,
the cohomology class of the transverse K\"ahler
metrics arising from elements in ${\mathcal S}(\xi, \bar{J})$  
is fixed, and it is natural to ask if there is a way of fixing the 
affine parameters $\varphi$ and $\psi$ also, which would yield then a 
canonical representative of this set. We proceed to discuss and answer this 
problem.

We start by introducing a Riemannian functional 
whose critical point will fix a canonical choice of metric for structures 
in ${\mathcal S}(\xi, \bar{J})$. This, in effect, will fix the desire
preferred representative that we seek.

We denote by ${\mathfrak M}(\xi,\bar{J})$ the set of all compatible 
Sasakian metrics arising from structures in ${\mathcal S}(\xi, \bar{J})$, and 
define the functional
\begin{equation}
\begin{array}{ccl}
{\mathfrak M}(\xi,\bar{J})  & \stackrel{E}{\rightarrow} & {\mathbb R}\, , \\
g & \mapsto & {\displaystyle \int _M s_g ^2 d{\mu}_g }\, , \end{array}
\label{var}
\end{equation}
the squared $L^2$-norm of the scalar curvature of $g$.

The variation of a metric in ${\mathfrak M}(\xi,\bar{J})$ depends 
upon the two {\it affine
parameters} of freedom $\varphi$ and $\psi$ of Lemma \ref{le1}. However,
the transversal K\"ahler metric varies as a function of $\varphi$ only, and
does so within a fixed basic cohomology class. The critical point of 
(\ref{var}), should it exist, will allow us to fix the parameter $\varphi$ 
since, not surprisingly, it shall be determined by the condition that 
$d\tilde{\eta}=d\eta +i\ddb \v$ be an extremal K\"ahler metric \cite{cal1}
on ${\mathcal D}$.
The remaining gauge function parameter $\psi$ represents nothing more than a 
change of coordinates in the representation of the form $\tilde{\eta}$ of 
the Sasakian structure in question. Thus, the finding of a 
critical point of $E$ produces a canonical representative of 
${\mathcal S}(\xi, \bar{J})$.

\subsection{Variational formulae}
In order to derive the Euler-Lagrange equation of (\ref{var}), we describe
the infinitesimal variations of the volume form, Ricci tensor, and scalar 
curvature, as a metric in ${\mathfrak M}(\xi,\bar{J})$ is deformed within that
space.

We begin by recalling that the Ricci form $\rho$ of a Sasaki 
structure $(\xi,\eta,\Phi,g)$ is defined on ${\mathcal D}$ by the 
expression
$$\rho_g(X,Y)=Ric_g(JX,Y) \, .$$
This is extended trivially on the characteristic foliation
$L_{\xi}$, and by Proposition \ref{pr1}, we easily see that 
\begin{equation}
\rho_g = \rho_g^T -2d\eta \, ,\label{rf}
\end{equation}
where $\rho^T$ denotes the form arising from the Ricci tensor of the
transversal metric, a basic two form that we think of as a 
$\bar{J}$-invariant two tensor in $\nu({\mathcal F}_{\xi})$. Thus, 
$\rho_g$ induces a well defined bilinear map on 
$\nu({\mathcal F}_{\xi})$ that is $\bar{J}$-invariant.

Though the notation suggests so, it is not the case that the
trace of the form (\ref{rf}) yields the scalar curvature of $g$. This
form only encodes information concerning second covariant derivatives of
$g$ along directions in $\nu({\mathcal F}_{\xi})$.

\begin{proposition} \label{pr2}
Let $(\xi,\eta_t,\Phi_t,g_t)$ be a path in ${\mathcal S}(\xi,\bar{J})$ that 
starts at $(\xi,\eta,\Phi,g)$ when $t=0$, and is such that
$d\eta_t=d\eta+ti\ddb \varphi$ for certain basic function $\v$, and for 
$t$ sufficiently small. Then we have the expansions
$$
\begin{array}{rcl}
d\mu_t & = & (1-\frac{t}{2}\Delta_B \varphi ))d\mu +O(t^2) \, , \\
\rho_t & = & \rho -ti\ddb \left( \frac{1}{2}\Delta_B \v + \v\right) + 
O(t^2)\, , \\
s_t &  = & s^T-2n -t\left( \frac{1}{2}\Delta_B^2 \varphi+2(\rho^T,
i\ddb \varphi)\right) +O(t^2) \, ,
\end{array}
$$ for the volume form, Ricci form, and scalar curvature of $g_t$,
respectively. Here, the geometric terms without sub-index are those 
corresponding to the starting metric $g$, and $\Delta_B$ is the Laplacian
acting on basic functions. 
\end{proposition}

{\it Proof}. By Lemma \ref{le1}, there exists a function $\psi$ such that
$$\eta_t=\eta+t(d^c \varphi + d\psi) + O(t^2) \, .$$
Since $\varphi$ and $\psi$ are both basic, we have that
$d\mu_t = \frac{1}{n!}\eta_t \wedge (d\eta_t)^n=\eta \wedge (d\eta_t)^n$,
and we obtain 
$$d\mu_t = \frac{1}{n!}\eta\wedge (d\eta +ti\ddb \varphi)^n=
d\mu +\frac{t}{(n-1)!}\eta \wedge
 (d\eta)^{n-1} \wedge i\ddb \varphi +O(t^2)\, .$$
Now, $\o_B=d\eta$, the K\"ahler form of the induced metric on 
${\mathcal D}$, is basic. We then have that $*_B (\o_B)^{n-1}/(n-1)!=
\o_B$, and conclude that
$$d\mu_t =\left( 1 -\frac{t}{2}\Delta_B \varphi \right) d\mu +O(t^2)\, .$$
 
By (\ref{rf}), we may compute the variation of $\rho$ by computing the
variation of $\rho^T$. This is well known to be \cite{si}
$$\rho_t^T=\rho^T +\frac{t}{2}i\ddb ( \Delta_B \varphi ) + O(t^2)\, .$$
Since $d\eta_t=d\eta+ti\ddb \v$, we obtain
$$\rho_t = \rho + ti\ddb \left( \frac{1}{2}\Delta_B \v + \v\right) + 
O(t^2) \, ,$$
as stated.

Finally, by Proposition \ref{pr1} once again, we have that the variation
of the scalar curvature arises purely from the variation of its transversal
part. Since the transversal metric is K\"ahler, we obtain
$$s_t=s-t \left( \frac{1}{2}\Delta_B^2 \varphi+2(\rho^T,i\ddb \varphi)\right)
+O(t^2) \, .$$
as desired. \qed

\begin{remark}
The forms $\rho^T$ and $i\ddb \v$ are basic. Hence, the metric pairing 
of these forms that appears in the Proposition above involves only the
transversal metric. On the other hand, in view of the analogous variational
formulae in the K\"ahler case \cite{si}, we might think that the
expression for $\rho_t$ above is a bit strange. We see that this is not so if
we just keep in mind that the Sasakian $\rho_g$ encodes 
$\nu({\mathcal F}_{\xi})$-covariant derivatives information only.  
\end{remark}

\subsection{Euler-Lagrange equations}\label{s2.2} 
Associated to any Sasakian structure $(\xi,\eta,\Phi,g)$ in 
${\mathcal S}(\xi,\bar{J})$, we introduce a basic 
differential operator $L_g^B$, of order $4$, whose kernel consists of 
basic functions with transverse holomorphic gradient. 

Given a {\it basic} function $\v: M\mapsto {\mathbb C}$,
we consider the vector field $\d^{\#}\v$ defined by the identity
\begin{equation}
g(\d^{\#}\varphi , \, \cdot \, ) = \db \varphi \, .\label{sh}
\end{equation}
Thus, we obtain the (1,0) component of the gradient of 
$\v$, a vector field that, generally speaking, is not transversally
 holomorphic.
In order to ensure that, we would need to impose
the condition $\bar{\partial}\partial^{\#} \v=0$, that is
equivalent to the fourth-order equation
\begin{equation}
(\bar{\d}{\d}^{\#})^{\ast}\bar{\d}{\d}^{\#}\v=0\, , \label{lic}
\end{equation}
because $\langle \v, (\bar{\d}{\d}^{\#})^{\ast}\bar{\d}{\d}^{\#} \v
\rangle_{L^{2}}=
\| \bar{\d}{\d}^{\#} \v\|^2_{L^{2}}$.

We have that
\begin{equation}
L_g^B\v := (\bar{\d}{\d}^{\#})^{\ast}\bar{\d}{\d}^{\#} \v 
=\frac{1}{4}(\Delta_B^2 \v + 4(\rho^T, i\ddb \v)+ 2({\partial}s^T )\hok
{\partial}^{\#}\v ) \, .\label{li}
\end{equation}
The functions on $M$ that are transversally constant are always in the
kernel of $L_g^B$. The only functions of this type that are basic are 
the constants. Thus, the kernel of $L_g^B$ has dimension at least $1$.

\begin{proposition} \label{prc}
The first derivative of $E$ at $g \in {\mathfrak M}(\xi,\bar{J})$ in the
direction of the deformation defined by $(\v, \psi)$ is given by
$$\frac{d}{dt}E(g_t)\mid_{t=0}=-4\int_{M}(s^T-2n) 
\left( (\bar{\d}{\d}^{\#})^{\ast}\bar{\d}{\d}^{\#}\v \right) \, d\mu \, .$$
\end{proposition}

{\it Proof}. This result follows readily from the fact that $s=s^T-2n$,
the variational formulae of Proposition \ref{pr2}, and identity
 (\ref{li}). \qed

As a corollary to Proposition \ref{prc}, we have:

\begin{theorem}\label{th1}
A Sasakian metric $g\in {\mathfrak M}(\xi,\bar{J})$ is a critical point of
the energy functional $E$ of {\rm (\ref{var})} if the basic 
vector field
$\d_g^{\#}s_g^T=\d_g^{\#}s_g$ is transversally holomorphic.
\end{theorem}

We are thus led to our fundamental definition:

\begin{definition}
We say that $\cals=(\xi,\eta,\Phi,g)$ is a {\bf canonical representative} of 
${\mathcal S}(\xi,\bar{J})$ if the metric $g$ satisfies the condition of 
Theorem {\rm \ref{th1}}, that is to say, if, and only if, $g$ is 
transversally extremal.
\end{definition}

\section{Transformation groups of Sasakian structures}
\setcounter{theorem}{0} 
In this section we discuss some important transformation groups associated to
Sasakian structures, and their corresponding Lie algebras.

Let us begin with a CR structure $({\mathcal D},J)$ of Sasaki type on $M$. If 
$\eta$ is a contact form, we have the group
$\gC\go\gn(M,\cald)$ of contact diffeomorphims, that is to say, the subgroup 
of the diffeomorphisms group $\gD\gi\gf\gf(M)$ consisting of those elements
that leave the contact subbundle $\cald$ invariant:
\begin{equation}
\gC\go\gn(M,\cald)=\{\phi\in \gD\gi\gf\gf(M)\, | \; \phi^*\eta=f_\phi\eta\, ,
 \;  f_\phi\in C^\infty(M)^*\}\, .
\end{equation}
Here, $C^\infty(M)^{*}$ denotes the subset of nowhere vanishing functions in 
$C^\infty(M)$. The dimension of the group so defined is infinite.

We may also consider the subgroup $\gC\go\gn(M,\eta)$ of strict contact
transformations, whose elements are those $\phi\in 
\gC\go\gn(M,\cald)$ such that $f_\phi=1$:
$$\gC\go\gn(M,\eta)=\{\phi\in \gD\gi\gf\gf(M)\, |\; \phi^*\eta=\eta \}\, .$$
This subgroup is also infinite dimensional. 

The Lie algebras of these two groups are quite important. The first of these
is the Lie 
algebra of {\it infinitesimal contact transformations},
\begin{equation}\label{infcon}
\gc\go\gn(M,\cald)=\{X\in \CX(M)\; | \; \pounds_X\eta=a(X)\eta\, , \; 
a(X)\in C^\infty(M)\}\, ,
\end{equation}
while the second is the subalgebra of {\it infinitesimal strict contact 
transformations}
\begin{equation}\label{strictcontacttrans}
\gc\go\gn(M,\eta)=\{X\in \CX(M)\; | \; \pounds_X\eta =0\}\, .
\end{equation}

If we now look at the pair $(\cald,J)$, we have the group of CR automorphisms
of $(\cald,J)$, defined by
\begin{equation}\label{crtrans}
\gC\gR(M,\cald,J)=\{\phi\in \gC\go\gn(M,\cald)\; | \;  \phi_*J=J\phi_*\}\, .
\end{equation}
This is a Lie group \cite{ChMo74}. Its Lie algebra $\gc\gr(M,\cald,J)$ 
can be characterized as
\begin{equation}\label{infcrtrans}
\gc\gr(M,\cald,J)=\{X\in \gc\go\gn(M,\cald)~|~ \pounds_XJ=0\}\, .
\end{equation}
Notice that in this defining expression, $\pounds_X J$ makes sense even though 
$J$ is not a tensor field on $M$, the reason being that the vector 
field $X$ leaves $\cald$ invariant.

Let us now consider a Sasakian structure $\cals=(\xi,\eta,\Phi,g)$ with 
underlying CR structure $({\mathcal D},J)$. We are interested in the subgroup 
of $\gD\gi\gf\gf(M)$ that leaves the tensor field $\Phi$ invariant. So we 
define
\begin{equation}
{\mathfrak S}_{\Phi}=\{\phi\in \gD\gi\gf\gf(M): \; \phi_*\circ\Phi=\Phi\circ 
\phi_*\}\, . \label{Phidiff}
\end{equation}
We also have
$$\gF\go\gl(M,\calf_\xi)=\{\phi\in \gD\gi\gf\gf(M): \; \phi_*\calf_\xi\subset 
\calf_\xi \}\, ,$$
the subgroup of $\gD\gi\gf\gf(M)$ that preserves the characteristic 
foliation of the Sasakian structure $\cals$.

In order to simplify the notation, we will often drop $M$ from the notation 
when referring to these various groups and algebras. 

\begin{lemma}\label{Philemma}
Let $(\cald,J)$ be a strictly pseudo-convex {\rm CR} structure of Sasaki type 
on $M$, and fix a Sasakian structure $\cals=(\xi,\eta,\Phi,g)$ with underlying 
CR structure $(\cald,J)$. Then
$$\gS_\Phi = \gC\gR(\cald,J)\cap \gF\go\gl(\calf_\xi)\, .$$
\end{lemma}

{\it Proof}. If $\phi\in\gS_\Phi$, the identity $\eta \circ \Phi =0$ implies 
that $\phi$ preserves $\cald$. If $X$ is a section of $\cald$, then we have 
that $\phi_*J(X)=\phi_*\Phi(X)=\Phi(\phi_*X)=J(\phi_*X)$, which implies that
$\phi\in \gC\gR(\cald,J)$.  But we have $\phi_*\Phi(\xi)=0=\Phi(\phi_*\xi)$
also, which implies that $\phi\in \gF\go\gl(\calf_\xi)$.

Conversely, suppose that $\phi\in \gC\gR(\cald,J)\cap \gF\go\gl(\calf_\xi)$.
Then $\phi$ leaves all three $\cald$, $\calf_\xi$, and $J$ invariant, and
therefore, it preserves the splitting (\ref{cs}).
Since the relation between $J$ and $\Phi$ is given by
\begin{equation} \label{Phidef}
\Phi(X)=\left\{ \begin{array}{l}
J(X) \; \text{if $X$ is a section of $\cald$}\, , \\
 0\; \text{if $X=\xi$}\, ,
\end{array} \right.
\end{equation} 
in order to conclude that $\phi \in \gS_\Phi$, it suffices to show that 
$\phi_*\Phi(\xi)=\Phi(\phi_*\xi)$. But this is clear as both sides 
vanish. \qed

Since $\gS_\Phi$ is closed in the Lie group $\gC\gR(\cald,J)$, it is itself
a Lie group. Generally speaking, the inclusion $\gS_{\Phi}\subset 
\gC\gR({\mathcal D}, J)$ is strict, and the group $\gF\go\gl(M,\calf_\xi)$ is 
infinite dimensional. 

The automorphism group $\gA\gu\gt(\cals)$ of the Sasakian structure
$\cals=(\xi,\eta,\Phi,g)$ is defined to be the subgroup of $\gD\gi\gf\gf(M)$ 
that leaves all the tensor fields in $(\xi,\eta,\Phi,g)$ invariant. 
It is a Lie group, and one has natural group inclusions
\begin{equation}\label{autcr}
\gA\gu\gt(\cals)\subset  \gS_\Phi\subset \gC\gR({\mathcal D}, J)
\end{equation}
whenever the CR structure $(\cald,J)$ is of Sasaki type, and $\cals$ has it
as its underlying CR structure.

The Lie algebras of $\gS_\Phi$ and $\gF\go\gl(\calf_\xi)$ are given by
\begin{equation}
{\mathfrak s}_{\Phi}=\{X\in \CX(M)\, | \; \pounds_X \Phi=0\}\, , \label{PhiLie}
\end{equation}
and
\begin{equation}\label{foliate}
\gf\go\gl(\calf_\xi) =\{X\in \CX(M)\, |\; [X,\xi]\; \text{is tangent to 
the leaves of $\calf_\xi$}\}\, ,
\end{equation}
respectively. The latter is just the Lie algebra of {\it foliate} vector 
fields of the foliation $\calf_\xi$. On the other hand, we may 
restate the defining condition for ${\mathfrak s}_{\Phi}$ as
$$
X\in {\mathfrak s}_{\Phi}  \Longleftrightarrow
 [X,\Phi(Y)]=\Phi([X,Y]) \;
\text{\rm for all} \; Y\, , 
$$
and we can see easily that
\begin{equation}\label{crphi}
\gs_\Phi = \gc\gr(\cald,J)\cap \gf\go\gl(\calf_\xi)\, .
\end{equation}

We can characterize now CR structures of Sasaki type in terms of their 
relations to the Lie algebras above.

\begin{lemma}\label{Reebcr}
Let $(\cald,J)$ be a strictly pseudoconvex {\rm CR} structure on $M$, and 
let $\eta$ be a compatible contact form representing $\cald$, with Reeb vector 
field $\xi$. Define the tensor field $\Phi$ by Equation {\rm (\ref{Phidef})}.
Then $(\cald,J)$ is of Sasaki type if, and only if, $\xi\in 
\gc\gr(\cald,J)$.
\end{lemma}

{\it Proof}. Given a strictly pseudoconvex CR structure $(\cald,J)$,
with contact 1-form $\eta$ that represents $\cald$, and Reeb 
vector field $\xi$, we consider the $(1,1)$ tensor field $\Phi$ defined by
(\ref{Phidef}). By Proposition 3.5 of \cite{BG01a}, $(\xi,\eta,\Phi)$ defines 
a Sasakian structure if, and only if, the CR structure is integrable and 
$\pounds_\xi\Phi =0$. Since $\xi$ is a foliate vector field, the condition 
that $\xi\in \gc\gr(\cald,J)$ ensures that $\xi \in {\mathfrak s}_{\Phi}$, 
and the result follows. \qed

For any strictly pseudoconvex CR structure with contact 1-form $\eta$, the 
Reeb vector field $\xi$ belongs to $\gc\go\gn(M,\eta)$. Therefore, if 
$\xi\in \gs_\Phi$, the structure $\cals=(\xi,\eta,\Phi,g)$ is Sasakian,
and $\xi\in \ga\gu\gt(\cals)$.

The study of the group $\gC\gR(\cald,J)$ has a long and ample history 
\cite{web77,Lee96,Sch95}. We state the most general result given in 
\cite{Sch95}. This holds for a general CR manifold $M$, so we
emphasize its consequence in the case where $M$ is closed.

\begin{theorem}\label{schthm}
Let $M$ be a $2n+1$ dimensional manifold with a strictly pseudoconvex {\rm CR}
structure $(\cald,J)$, and {\rm CR} automorphism group $G$. If $G$ does not 
act properly on $M$, then:
\begin{enumerate}
\item If $M$ is a non-compact manifold, then it is {\rm CR} diffeomorphic to 
the Heisenberg group with its standard {\rm CR} structure.
\item If $M$ is a compact manifold, then it is {\rm CR} diffeomorphic to the 
sphere ${\mathbb S}^{2n+1}$ with its standard {\rm CR} structure.
\end{enumerate}
In particular, if $M$ is a closed manifold not {\rm CR} diffeomorphic to 
the sphere,
the automorphisms group of its {\rm CR} structure is compact. 
\end{theorem}

We recall that the inclusion $\gA\gu\gt(\cals)\subset \gC\gR({\mathcal D}, J)$
(see (\ref{autcr})), generally speaking, is proper. However, we 
do have the following.

\begin{proposition}\label{autSas}
Let $(\cald,J)$ be a strictly pseudoconvex {\rm CR} structure on a closed 
manifold
$M$, and suppose that $(\cald,J)$ is of Sasaki type. Then 
there exists a Sasakian structure $\cals=(\xi,\eta,\Phi,g)$ with underlying 
{\rm CR} structure $(\cald,J)$, whose automorphism group $\gA\gu\gt
(\cals)$ is a maximal compact subgroup of $\gC\gR(M,\cald,J)$. In fact,
except for the case when $(M,\cald,J)$ is {\rm CR} diffeomorphic to the sphere 
${\mathbb S}^{2n+1}$ with its standard {\rm CR} structure, the automorphisms 
group $\gA\gu\gt(\cals)$ of $\cals$ is equal to
$\gC\gR(M,\cald,J)$.
\end{proposition} 

{\it Proof}. Let $G$ be a maximal compact subgroup of $\gC\gR(\cald,J)$.
By Theorem \ref{schthm}, we have that $G=\gC\gR(\cald,J)$ except when 
$(M,{\mathcal D},J)$ is CR diffeomorphic to the sphere 
${\mathbb S}^{2n+1}$.
Let $\tilde{\cals}=(\tilde{\xi},\tilde{\eta},\tilde{\Phi},\tilde{g})$ be a 
Sasakian structure with underlying CR structure $(\cald,J)$. 

If $\phi\in G$, then $\phi^*\tilde{\eta}=f\tilde{\eta}$ for some nowhere 
vanishing real-valued function $f$. By averaging $\tilde{\eta}$ over $G$, we 
obtain a 
$G$-invariant contact form $\eta$ with associated contact structure ${\cald}$. 
Let $\xi$ be its Reeb vector field. As $\xi$ is uniquely determined by $\eta$, 
we conclude that $\phi_*\xi=\xi$. We then define a $(1,1)$-tensor 
$\Phi$ by the expression in (\ref{Phidef}). The conditions 
$\phi_*J=J\phi_*$ and $\phi_*\xi=\xi$ imply that $\phi_*\Phi=\Phi\phi_*$.
The triple $(\xi,\eta,\Phi)$ defines a Sasakian structure $\cals$, and
$\phi\in  \gA\gu\gt(\cals)$. \qed 
\medskip

We now look at the case where the manifold $M$ is polarized by
$(\xi, \bar{J})$. Then, $(\nu({\mathcal F}_{\xi}),\bar{J})$ is a complex
vector bundle, and any $\phi\in \gF\go\gl(M,\calf_\xi)$ induces a map 
$\bar{\phi}_*:\nu(\calf_\xi)\rightarrow \nu(\calf_\xi)$. We define the group 
of {\it 
transversely holomorphic transformations} $\gH^T(\xi,\bar{J})$ by 
\begin{equation}\label{transhologrp}
\gH^T(\xi,\bar{J})=\{\phi\in \gF\go\gl(M,\calf_\xi) ~|~ \bar{\phi}_*
\circ\bar{J}=\bar{J}\circ 
\bar{\phi}_*\}.
\end{equation}
Since a 1-parameter subgroup of any smooth section of $L_\xi$ induces the 
identity on 
$\nu(\calf_\xi)$, this group is infinite dimensional. We are mainly interested
 in the 
infinitesimal version. Given a choice of $\cals=(\xi,\eta,\Phi,g)$ in
${\mathcal S}(\xi,\bar{J})$, $(\nu({\mathcal F}_{\xi}),\bar{J})$ is identified 
with the underlying CR structure $({\mathcal D},J)$. In this case, we may use 
the decomposition (\ref{cs}) to write any vector field as
\begin{equation}
X=X_{\mathcal D}+c(X)\xi \, ,\label{dec}
\end{equation}
which defines the component function $X \mapsto c(X):=\eta(X)$, and the class
$\bar{X}$ on $\nu({\mathcal F}_{\xi})$ defined by the vector field $X$ is
represented by $X_{\mathcal D}$. 

The Lie bracket operation induces a bilinear mapping 
on $\nu({\mathcal F}_{\xi})$ by
$$[\bar{X},\bar{Y}]:=\overline{[X,Y]}\, .$$
This operation allows us to generalize the notion of transversally 
holomorphic vector field already encountered in \S\ref{s2.2}.

\begin{definition}
Let $(M,\xi,\bar{J})$ be a polarized Sasakian manifold. We say that a vector 
field $X$ is {\bf transversally  holomorphic} if given any section $\bar{Y}$ of
$\nu({\mathcal F}_{\xi})$, we have that 
$$[\bar{X},\bar{J}\, \bar{Y}]=\bar{J}\, \overline{[X,Y]}\, .$$
The set of all such vector fields will be denoted by 
${\mathfrak h}^T(\xi, \bar{J})$.
\end{definition}

It is now desirable to express the defining condition for $X$ to be
in ${\mathfrak h}^T(\xi, \bar{J})$ in terms intrinsic to $X$ itself.
The reader may consult \cite{brsl} for relevant discussions.  

\begin{lemma}
Let $(M,\xi, \bar{J})$ be a polarized Sasakian manifold. If 
$X\in {\mathfrak h}^T(\xi, \bar{J})$ is a transversally holomorphic
vector field, for any Sasakian structure $(\xi,\eta,\Phi,g)\in
{\mathcal S}(\xi,\bar{J})$ we have that
$$(\pounds_X \Phi)(Y)=\eta([X,\Phi(Y)])\xi \, .$$
\end{lemma}

{\it Proof}. We have that 
$$
(\pounds_X \Phi)(Y)=[X,\Phi(Y)]-\Phi([X,Y])\, ,$$
which implies that $c((\pounds_X \Phi)(Y))=\eta([X,\Phi(Y)])$. Thus,
$$(\pounds_X \Phi)(Y)=((\pounds_X \Phi)(Y))_{\mathcal D}+
\eta([X,\Phi(Y)])\xi\, .$$
The result follows after simple considerations.  \qed

The set $\gh^{T}(\xi,\bar{J})$ is a Lie algebra contained in
$\gf\go\gl(M,\calf_\xi)$. If we represent $(\nu({\mathcal F}_{\xi}),\bar{J})$ 
as  $({\mathcal D},J)$ for a choice of $\cals=(\xi,\eta,\Phi,g)$ in
${\mathcal S}(\xi,\bar{J})$ with underlying CR structure $({\mathcal D},J)$,
by the decomposition (\ref{dec}) we see that for a
transversally holomorphic vector field $X$ we have that
$$[X_{\mathcal D},J(Y_{\mathcal D})]_{\mathcal D}=
J([X_{\mathcal D},Y_{\mathcal D}]_{\mathcal D}) \, ,$$
for any vector field $Y$. Thus, $X_{\mathcal D}$ preserves the transverse 
complex structure $J$. This characterization can be reformulated by 
saying that if $X\in {\mathfrak h}^T(\xi,\bar{J})$, the vector field of 
type $(1,0)$ given by
\begin{equation}
\Xi_X=\frac{1}{2}(X_{\mathcal D}-iJ(X_{\mathcal D})) \label{hc}
\end{equation}
is in the kernel of the transverse Cauchy-Riemann equations. Thus, the 
mapping into the space of sections of $\nu({\mathcal F}_{\xi})$ given
by  
\begin{equation}
\begin{array}{ccc}
\gh^T(\xi,\bar{J}) & \rightarrow & \Gamma (\nu({\mathcal F}_{\xi})) \\
X & \mapsto & \bar{X} 
\end{array}
\end{equation}
has an image that can be identified with the space of sections of 
$(\cald,J)$ satisfying the Cauchy Riemann equations, there finite dimensional.
We denote this image by $\gh^T(\xi,\calf_\xi)/L_\xi$. 

Their one dimensional foliations made Sasakian manifolds a bit special. In
particular, they carry no non-trivial parallel vector field (see \cite{BG06}).
For by a result of Tachibana, any harmonic one form must annihilate
the Reeb vector field $\xi$, and so any parallel vector fields $X$ must be 
orthogonal to $\xi$, that is to say, it must be a section of ${\mathcal D}$. 
But then, since the metric is covariantly constant and $\xi$ is Killing, we
must have that $0=\nabla_Y g(\xi,X)=g(\Phi(Y),X)$ for all $Y$, which 
forces $X$ to be identically zero.

\begin{remark}\label{MatLic}
Notice that for any Sasakian structure with underlying CR structure
$(\cald,J)$, the Lie algebra $\gc\gr(\cald,\eta)$ is reductive. For
either $\gC\gR(\cald,J)$ is compact, or $\gc\gr(\cald,\eta)=\gs\gu(n+1,1)$.
In particular, let $(N,\gro)$ be a K\"ahler manifold, and consider the circle 
bundle $\pi:M \rightarrow N$ with Euler class $[\gro]$. Then $M$ has a natural
Sasakian structure $(\xi,\eta,\Phi,g)$ such that $\pi^*\gro=d\eta$. Then the
horizontal lifts of non-trivial parallel vector fields in $(N,\omega)$ are 
holomorphic but not parallel on $M$. On the other hand, if we take a 
holomorphic field that lies in the non-reductive part of the
algebra of holomorphic vector fields of $N$, its horizontal lift is
a transversally holomorphic vector field that does not lie in the reductive
component of the algebra ${\mathfrak h}^T$, and thus, it cannot
possibly leave the contact subbundle $\cald$ invariant. \qed
\end{remark}
\medskip

If $X\in \gh^T(\xi,\bar{J})$, given any real valued function $f$,
$X+f\xi\in \gh^T(\xi,\bar{J})$ also, and so, $\gh^T(\xi,\bar{J})$ cannot have
finite dimension. The remark above alludes to the special structure that 
$\gh^T(\xi,\bar{J})/L_{\xi}$ has, and in fact, we are now
ready to extend to the Sasakian context a result of Calabi \cite{cal2} on 
the structure of the algebra of holomorphic vector fields of a K\"ahler 
manifold that carries an 
extremal metric. Calabi's theorem is, in turn, an extension of work of 
Lichnerowicz \cite{Lic57} on constant scalar curvature metrics, and the latter
is itself an extension of a result of Matsushima \cite{Mat57} in the 
K\"ahler-Einstein case. We also point the reader to the theorem for harmonic 
K\"ahler foliations in \cite{nito}, which is relevant in this context.

Consider a Sasakian structure $(\xi,\eta,\Phi,g)$ in 
${\mathcal S}(\xi,\bar{J})$.
Let ${\mathcal H}^B_g$ be the space of basic functions in the kernel of the 
operator $L_g^B$ in (\ref{li}), and consider the mapping
\begin{equation}
\partial_{g}^{\#}: {\mathcal H}_g^B \rightarrow \gh^T(\xi,
\bar{J})/L_{\xi}
\, , \label{ma1}
\end{equation}
where $\partial_g^{\#}$ is the operator defined in (\ref{sh}).
We use the Sasakian metric $g$ to identify the quotient space in the right 
side above with the holomorphic vector fields that are sections of 
$({\mathcal D},J)$, which we shall refer
to from here on as ${\mathfrak h}({\xi,\mathcal D},J)$. The notation for 
this Lie algebra is a bit 
non-standard in that ${\mathfrak h}({\xi,\mathcal D},J)$ depends on 
$\xi$ or rather on the foliation 
$\calf_\xi$ and are holomorphic sections of $\cald$ with respect to $J$.
It should be noted, however, that 
while elements in ${\mathfrak h}({\xi,\mathcal D},J)$ leave both 
$\calf_\xi$ and $\bar{J}$ invariant, they do not necessarily leave 
$\cald$ invariant. 

We also define the 
operator $\bar{L}^B_g$ on ${\mathcal H}^B_g$ by $\bar{L}^B_g \varphi =
\overline{L^B_g \bar{\varphi}}$. It follows that
\begin{equation}
(\bar{L}^B_g-L^B_g)\varphi = \partial_{g}^{\#}s_g \hok \partial \varphi-
\partial_{g}^{\#}\varphi \hok \partial s_g \, , \label{ma2}
\end{equation}
where $s_g$ is the scalar curvature of $g$. The fact that $s_g$ is a basic 
function implies that $\partial_{g}^{\#}s_g$ is a $(1,0$) section of
$({\mathcal D},J)$. The identity above implies that $L^B$ and $\bar{L}^B$ 
coincide if $s_g$ is constant. If the metric $g$ is canonical,
then we have that $\partial_{g}^{\#}s_g \in {\mathfrak h}(\xi,{\mathcal D},J)$,
and the operators $L^B$ and $\bar{L}^B$ commute.

The image ${\mathfrak h}_0\cong {\mathcal H}_g^B /{\mathbb C}$ of the mapping
(\ref{ma1}) is an ideal in ${\mathfrak h}(\xi,{\mathcal D},J)$, and can be
identified with the space of holomorphic fields that have non-empty zero
set. The quotient algebra ${\mathfrak h}(\xi,{\mathcal D},J)/{\mathfrak h}_0$ 
is  Abelian. We also denote by $\ga\gu\gt(\bar{J},g_T)$ the Lie subalgebra 
of $\gh^T(\xi,\nu(\calf_\xi))$ 
that are holomorphic Killing vector fields of the transverse metric 
$g_T$, that is
\begin{equation}\label{transKilling}
\ga\gu\gt(\bar{J},g_T)=\{\bar{X}\in \gh(\xi,\cald,J) ~|~ \pounds_{\bar{X}}
g_T=0 \}\, .
\end{equation}

Suppose now that $(\xi,\eta,\Phi,g)$ is a canonical representative of
${\mathcal S}(\xi,\bar{J})$, so that $g$ is Sasaki extremal. Let 
${\mathfrak z}_0$ be the image under $\partial_{g}^{\#}$ of the set of purely 
imaginary functions in ${\mathcal H}_g^B$. This is just the space of
Killing fields for the transversal metric $g^T$ that are of the form 
$J\nabla_{g^T} \varphi$, $\varphi \in {\mathcal H}^B_g$.
Furthermore, by (\ref{ma2}) we see
that the complexification ${\mathfrak z}_0\oplus \bar{J}{\mathfrak z}_0$ 
coincides with the commutator of $\partial_{g}^{\#}s_g$:
$${\mathfrak z}_0\oplus \bar{J}{\mathfrak z}_0=
\{ X \in {\mathfrak h}(\xi,{\mathcal D},J)\, : \;
 [X,\partial_{g}^{\#}s_g]=0\}\, .$$

\begin{theorem}\label{ml}
Let $(M,\xi,\bar{J})$ be a polarized Sasakian manifold. Suppose that
there exists a canonical representative $(\xi,\eta,\Phi,g)$ of
${\mathcal S}(\xi,\bar{J})$. Let ${\mathcal H}^B_g$ be the space of basic
functions in the kernel of the operator $L_g^B$ in {\rm (\ref{li})}, and let
${\mathfrak h}_0$ be the image of the mapping {\rm (\ref{ma1})}. Then
we have the orthogonal decomposition
$$\gh^T(\xi,\bar{J})/L_{\xi}\cong {\mathfrak h}(\xi,{\mathcal D},J)
={\mathfrak a}\oplus {\mathfrak h}_0\, ,$$
where ${\mathfrak a}$ is the algebra of parallel vector fields of the
transversal metric $g^T$. Furthermore,
$${\mathfrak h}_0={\mathfrak z}_0\oplus \bar{J}{\mathfrak z}_0\oplus 
(\oplus_{\lambda >0}{\mathfrak h}^{\lambda})\, ,$$
where ${\mathfrak z}_0$ is the image of the purely imaginary elements of
${\mathcal H}^B_g$ under $\partial_{g}^{\#}$, 
and ${\mathfrak h}^{\lambda}=\{ \bar{X} \in \gh^T(\xi,\bar{J})/L_{\xi}\, :
\; [\bar{X},\partial_{g}^{\#}s_g]=\lambda \bar{X}\}$.
Moreover, $\gz_0$ is isomorphic to the quotient algebra 
$\ga\gu\gt(\xi,\eta,\Phi,g)/\{\xi\},$ 
so the Lie 
algebra $\ga\gu\gt(\bar{J},g_T)$ of Killing vector fields
for the transversal metric $g^T$ is equal to
$$\ga\gu\gt(\bar{J},g_T)={\mathfrak a}\oplus {\mathfrak z}_0\cong {\mathfrak 
a}\oplus\ga\gu\gt(\xi,\eta,\Phi,g)/\{\xi\} \, .$$
\end{theorem}

The presence of the algebra ${\mathfrak a}$ above does not contradict the 
fact that there are no non-trivial parallel vector fields on a closed Sasakian 
manifold: a vector field can be parallel with respect to $g^T$ without being 
parallel with respect to $g$.

{\it Proof of Theorem} \ref{ml}. We prove the last statement first. In order 
to see this, we notice that there is an 
exact sequence \cite{BG06}
\begin{equation}\label{obstructseq}
0 \rightarrow \{\xi\} \rightarrow \ga\gu\gt(\xi,\eta,\Phi,g) \rightarrow 
\ga\gu\gt(\bar{J},g_T) \stackrel{\grd}{\rightarrow}
H^1_B(\calf_\xi)\, ,
\end{equation}
where $H^1_B(\calf_\xi)$ denotes the basic cohomology associated to the 
characteristic foliation 
$\calf_\xi$. Using the identification of $\ga\gu\gt(\bar{J},g_T)$ with 
elements in 
$\gh(\xi,\cald,J),$ let us 
describe the map $\grd$. Since $\bar{X}\in \ga\gu\gt(\bar{J},g_T)$ it leaves 
$d\eta$ invariant, so the 
1-form $\bar{X}\hok d\eta$ is closed and basic. It, thus, defines an element 
in 
$H^1_B(\calf_\xi)$. 
So we can define $\grd(\bar{X})=[\bar{X}\hok d\eta]_B$. Now the section 
$\bar{X}\in 
\ga\gu\gt(\bar{J},g_T)$ can be extended to an element $X=\bar{X}+a\xi\in 
\ga\gu\gt(\xi,\eta,\Phi,g)$ 
if, and only if, the basic cohomology class $[\bar{X}\hok d\eta]_B$ vanishes, 
and this determines $a$ 
up to a constant. By Hodge theory and duality, the image of $\grd$ can be 
identified with the Lie algebra of 
parallel vector fields in $\ga\gu\gt(\bar{J},g_T)$. The splitting then follows 
as in the K\"ahler case \cite{cal2}.

For the first part of the theorem we sketch the main points, as the argument 
is an adaptation to
our situation of that in \cite{cal2}. Given a section $\bar{X}$ in 
${\mathfrak h}(\xi,{\mathcal D},J)$, we look at the Hodge decomposition of the
$(0,1)$-form that corresponds to it via the metric $g^T$. It is $\db$-closed,
and both, its harmonic and $\db$ components, are the dual of holomorphic
fields. The vector field dual to the harmonic component is $g^T$-parallel.
 
Since $g^T$ is an extremal metric, the operators $L_g^B$ and $\bar{L}_g^B$
commute. We then restrict $\bar{L}_g^B$ to the kernel of $L_g^B$, and
use the resulting eigenspace decomposition together with the identity
(\ref{ma2}) to derive the remaining portion of the theorem. \qed
 
\begin{remark}
This result obstructs the existence of special canonical representatives of 
a polarized Sasakian manifold in the same way it does so in the K\"ahlerian
case. For instance, let $(N,\gro)$ be the one-point or two-points blow-up
of ${\mathbb C}{\mathbb P}^2$, and consider the circle 
bundle $\pi:M \rightarrow N$ with Euler class $[\gro]$. If  $M$ is polarized
by its natural Sasakian structure $(\xi,\eta,\Phi,g)$, the one where
$\pi^*\gro=d\eta$, then ${\mathcal S}(\xi,\bar{J})$ cannot be represented
by a Sasakian structure $(\xi,\tilde{\eta},\tilde{\Phi},\tilde{g})$ with 
$\tilde{g}$ a metric of constant scalar curvature. The structure of 
${\mathfrak h}^T(\xi,\bar{J})/L_{\xi}$ would obstruct it. \qed
\end{remark}

\section{A Sasaki-Futaki invariant}
\setcounter{theorem}{0} 
Let $(M,\xi,\bar{J})$ be a polarized Sasakian manifold. Given any 
structure $(\xi,\eta,\Phi,g)\in {\mathcal S}(\xi,\bar{J})$, we denote
its underlying CR structure by $({\mathcal D},J)$. The metric $g$ is
an element of ${\mathfrak M}(\xi,\bar{J})$ whose transversal 
Ricci form $\rho^T$ is basic.
We define the Ricci potential $\psi_g$ as the function in the Hodge
decomposition of $\rho^T$ given by
$$\rho^T=\rho^T_h +i\ddb \psi_g \, ,$$
where $\rho^T_h$ is the harmonic representative of the foliated 
cohomology class represented by $\rho^T$. Notice that if $G_g^T$ is the
Green's operator of the transversal metric, we have that
$$\psi_g=-G_g^{T}(s_g^T)=-G_g^T(s_g^T-2n)=-G_g^T(s_g)=-G(s_g-s_{g,0})\, ,$$
where $s_g$ and $G_g$ are the scalar curvature and Green's operator of $g$, and
$s_{g,0}$ is the projection of $s_g$ onto the constants.
The sequence of equalities above follows by (c) of Proposition \ref{pr1},
which implies that $s_g=s_g^T-2n$ is a basic function. Thus, the Ricci
potential $\psi_g$ is itself a basic function.

On ${\mathfrak h}^T(\xi,\bar{J})$, we define the function
$$X\mapsto \int X(\psi_g) d\mu_g \, .$$
Since $\psi_g$ is basic, the
integrand in this expression can be fully written in terms of the
transversally holomorphic realization $\Xi_X$ (see (\ref{hc})) of $X$.

\begin{proposition} 
The mapping above only depends on the basic cohomology class 
represented by $d\eta$, and not on the particular transversal K\"ahler metric
induced by $g\in {\mathfrak M}(\xi,\bar{J})$ that is used to represent it.
\end{proposition}

{\it Proof}. We take a path $g_t$ in ${\mathfrak M}(\xi,\bar{J})$ starting at
$g$ for which the
transversal K\"ahler form is of the form 
$$d\eta_t=d\eta + ti\ddb \v \, ,$$
with the affine parameter $\v$ a basic function.
From the identity $\Delta_B \psi_g=s^T_{g,0}-s^T$, we see that 
the variation $\dot{\psi}_g$ of $\psi_g$ satisfies the relation
$$2(i\ddb \v, i\ddb \psi_g)_g+\Delta_B \dot{\psi}_g=-\dot{s^T}=
\frac{1}{2}\Delta_B^2 \v+2(\rho^T,i\ddb \v)_g \, .$$
Hence,
$$\dot{\psi}_g-\frac{1}{v}
\int \dot{\psi}_g d\mu_g =\frac{1}{2}\Delta_B \v +2G^T_g(\rho_h^T,i\ddb \v)_g
 \, ,$$
where $v$ is the volume of $M$ in the metric.

Since $\rho_h^T$ is harmonic, the last summand in the right side
can be written as
$-2G^T_g(\d^{*}(\db^{*}(\v \rho_h^T)))$. For convenience, let us set
$\beta=\db^{*}(\v \rho_h^T)$. Hence,
$$\frac{d}{dt}\int X(\psi_t) d\mu_{g_t} =
\int X\left( \frac{1}{2}\Delta_B \v - 2G^T_g(\d^{*}\beta))-\frac{1}{2}
\psi \Delta_B \v\right) d\mu_g \, .$$
By the Ricci identity for the transversal metric, we have that
$$\frac{1}{2}(\Delta_B \v)_{\a}=-\v_{,\g \a}^{\hspace{4mm}\g}
+ \v_{,\g}(\psi_{,\a}^{\hspace{2.1mm} \g}+(r_{h}^T)_{\a}^{\hspace{1mm} \g})=
-\v_{,\g \a}^{\hspace{4mm}\g}+ \v_{,\g}\psi_{,\a}^{\hspace{2.1mm} \g}+
\beta_{\a}
\, ,$$
and so, after minor simplifications, we conclude that
$$
\begin{array}{rcl}
{\displaystyle \frac{d}{dt}\int X(\psi_t) d\mu_{g_t}} & = & {\displaystyle 
\int \Xi_X^{\a}\left( \v_{, \g}
\psi_{\a}-\v_{,\g \a}\right)^{\g} d\mu_g \; +} \\ & & {\displaystyle 
\mbox{}\hspace{.5in} \int \Xi_X^{\a}\left( \beta_{\a} - 
2(G_g\d^{*}\beta)_{,\a}\right) d\mu_g }\, ,
\end{array}
$$
where $\Xi_X$ is the (1,0)-component of $X_{\mathcal D}$ (see (\ref{hc})).
 
The first summand on the right above is zero because $\Xi_X$ is holomorphic.
This is just a consequence of Stokes' theorem. The second summand is also
zero since we have
$$\begin{array}{rcl}
{\displaystyle \int \Xi_X^{\a}\left( \beta_{\a} - 2(G_g\d^{*}\beta)_{,\a}
\right) d\mu_g } & = & {\displaystyle 
\int (\beta-\Delta_B G^T_g \beta , \Xi_X^{\flat})d\mu_g }\; + \\ & &
 {\displaystyle \int (2\d^{*}\d G^T_g \beta, \Xi_X^{\flat})
d\mu_g }\, ,
\end{array}$$
and $\beta - \Delta_B G_g^T \eta = 0$ while $\d \Xi_X^{\flat }=0$. Here,
of course, $\Xi_X^{\flat }$ is the $(0,1)$-basic form corresponding the
$(1,0)$-vector field $\Xi_X$. \qed

We may then define the transversal Futaki invariant ${\mathfrak F}=
{\mathfrak F}_{(\xi,\bar{J})}$ of the polarized Sasakian 
manifold $(M,\xi,\bar{J})$ to be the functional
\begin{equation}
\begin{array}{c}
{\mathfrak F}: {\mathfrak h}^T(\xi,\bar{J}) 
\longrightarrow {\mathbb C} \\
{\mathfrak F} (X)={\displaystyle \int _{M}
X(\psi _{g})d\mu_g = -\int _{M} X(G_g^Ts_g^T) d\mu_g }\, ,\end{array}
\label{fusa} \end{equation}
where $g$ is any metric in ${\mathfrak M}(\xi,\bar{J})$.
The Proposition above shows that ${\mathfrak F}$ is well-defined, as this
expression depends only on the basic class $[d\eta]$ of a Sasakian structure
$(\xi,\eta, \Phi,g)$ in ${\mathcal S}(\xi,\bar{J})$, rather than the specific 
Sasakian structure chosen to represent it. 

It is rather obvious that ${\mathfrak F} (X)={\mathfrak F} (X_{\mathcal D})$,
and the usual argument in the K\"ahler case implies also that
${\mathfrak F} ([X,Y])=0$ for any pair of vector fields $X$, $Y$ in 
${\mathfrak h}^T$. 

Our next proposition extends to canonical Sasakian metrics a now
well-known result in K\"ahler geometry originally due to Futaki \cite{fu}.
The Sasaki version here is analogous to the expanded version of Futaki's result
presented by Calabi \cite{cal2}.

\begin{proposition}\label{FSscalar} 
Let $(\xi,\eta,\Phi,g)$ be a canonical Sasakian representative of 
${\mathcal S}(\xi, J)$. Then, the metric $g$ has constant scalar curvature 
if, and only if, ${\mathfrak F}(\, \cdot \, )=0$.
\end{proposition}

{\it Proof}. In one direction the statement is obvious: a constant
scalar curvature Sasakian metric has trivial Ricci potential function, 
and the functional (\ref{fusa}) vanishes on any $X\in {\mathfrak h}^T
(\xi,\bar{J})$.

In order to prove the converse, we first observe that if 
$X \in {\mathfrak h}^T(\xi,\bar{J})$ is a transversally holomorphic vector 
field of the form 
$X=\d_g^{\#}f$ for some basic function $f$, then 
$$\begin{array}{rcl}
{\displaystyle
{\mathfrak F} (X)=-\int \d_g^{\#}f (G_gs_g)d\mu_g } & = & {\displaystyle
-2\int (\db f, \db G^T_gs_g)_gd\mu_g } \\ & = & {\displaystyle 
-2\int f(\db_g^{*}\db G^T_gs_g) 
d\mu_g }\, ,
\end{array}$$
because the scalar curvature $s_g$ is a basic function also. Since 
$2\db^{*}\db= \Delta_B$, we conclude that
$${\mathfrak F} (\d_g^{\#}f)=- \int f(s_g-s_{g,0})d\mu_g \, .$$
Now, if the Sasakian metric $g$ is a critical point of 
the energy function $E$ in (\ref{var}), then $\d_g^{\#}s=\d^{\#}s^T=
\d_g^{\#}s_g$ is a transversally holomorphic vector 
field, and we conclude that 
$${\mathfrak F} (\d_g^{\#}s^T )=-\int (s_g-s_{g,0})^2 d\mu_g \, .$$
Thus, if ${\mathfrak F}(\, \cdot \, )=0$, then $s_g$ must be constant.
\qed

A particular case of constant scalar curvature Sasakian metrics is the case 
of Sasakian $\eta$-Einstein 
metrics. These are Sasakian metrics $g$ that satisfy
\begin{equation}\label{etaEin}
\Ric_g=\lambda g+\nu\eta\otimes \eta
\end{equation}
for some constants $\lambda$ and $\nu$. The scalar curvature 
$s_g$ of these metrics is given by
$s_g=2n(1+\lambda)$.
We refer the reader to \cite{BGM06}, and references therein, for further 
discussion of these type of metrics.  

\begin{corollary}
Let $(\xi,\eta,\Phi,g)$ be a canonical representative of 
${\mathcal S}(\xi, \bar{J})$, and suppose that the basic first Chern class 
$c_1({\mathcal F}_{\xi})$ is a 
constant multiple, say $a,$ of $[d\eta]_B$. Then
\begin{enumerate}
\item If $a=0$, then $(\xi,\eta,\Phi,g)$ is a null $\eta$-Einstein Sasakian 
structure with $\lambda=-2$, whose transverse metric is Calabi-Yau.
\item If $a<0$, then $(\xi,\eta,\Phi,g)$ is a negative $\eta$-Einstein 
Sasakian structure with $\lambda<-2$, whose transverse metric is 
K\"ahler-Einstein with negative scalar curvature.
\item If $a>0$, then $(\xi,\eta,\Phi,g)$ is a positive $\eta$-Einstein 
Sasakian structure with $\lambda>-2$, whose transverse metric is positive 
K\"ahler-Einstein if, and only if, the Futaki-Sasaki invariant 
${\mathfrak F}_{\xi,\bar{J}}$ vanishes. Moreover, if ${\mathfrak F}_{\xi,
\bar{J}}$ vanishes, 
$g\in {\mathfrak M}(\xi, \bar{J})$ is Sasaki-Einstein if, and only if, 
$\lambda=2n$.
\end{enumerate}

\noindent When an $\eta$-Einstein metric exists, the relation 
$2\pi a=\lambda+2$ holds.
\end{corollary}

{\it Proof}. Parts (1) and (2) follow from Proposition \ref{FSscalar} and 
Theorem 17 of \cite{BGM06}.  For 
(3) we notice that if $g$ is positive Sasakian $\eta$-Einstein, the result 
follows immediately from 
Proposition \ref{FSscalar}.
Conversely, if $g$ is a canonical representative ${\mathcal S}(\xi, \bar{J})$
and ${\mathfrak F}_{\xi,\bar{J}}$ vanishes, its scalar 
curvature is constant. It follows that the scalar curvature of the 
transversal metric is constant also, and this implies that 
$\rho_g+2d\eta = \rho^T$ 
is transversally harmonic. As the latter form represents 
$2\pi c_1({\mathcal F}_{\xi})$, which is also represented 
by a constant multiple of $d\eta$, the uniqueness of the harmonic 
representative of a class implies that $\rho_g +2d\eta= \rho^T=2\pi ad\eta$ 
for some
$a > 0$. It then follows from this that the transverse Ricci tensor $\Ric_T$ 
satisfies
$\Ric_T=2\pi ag_T$. But then $g$ is a positive Sasakian $\eta$-Einstein 
metric, and it follows from 
Equation \ref{etaEin} that $2\pi a=\lambda+2$. \qed

This Corollary applies whenever the first Chern class $c_1(\cald)$ of the 
contact bundle is a torsion class.
We mention also that one can always obtain a Sasaki-Einstein metric from a 
positive Sasakian 
$\eta$-Einstein metric by applying a transverse homothety.

\section{The Sasaki cone}\label{sc}
\setcounter{theorem}{0} 

We have the following result for CR structures of Sasaki type on $M$, 
an immediate consequence of the argument in the proof of Lemma \ref{Reebcr}.

\begin{proposition}\label{xicr}
Let $({\mathcal D},J)$ be a {\rm CR} structure of Sasaki type on $M$, and let 
$\cals=(\xi,\eta,\Phi,g)$  be a contact metric structure whose underlying 
{\rm CR} structure is $(\cald,J)$. Then $\cals$ is a Sasakian structure if
and only if $\xi\in \gc\gr(\cald,J)$.
\end{proposition}

We fix a strictly pseudoconvex CR structure $(\cald,J)$ on $M$, and
define the set 
\begin{equation}
{\mathcal S}(\cald,J)=\left\{
\begin{array}{c}
 \cals=(\xi,\eta,\Phi,
g):\; \cals \; {\rm a\; Sasakian\; structure} \\
({\rm ker}\, \eta, \Phi \mid_{{\rm ker}\, \eta})=
(\cald,J)\end{array} \right\}\, .
\end{equation}
We think of this as a subspace of sections of a vector bundle, and provide it
with the $C^\infty$ compact-open topology. This set
is nonempty if, and only if, $(\cald,J)$ is of Sasaki type.

\begin{proposition}\label{poseta}
Let $(\cald,J)$ be a {\rm CR} structure of Sasaki type, and let 
$\cals_0=(\xi_0,\eta_0,\Phi_0,g_0)\in {\mathcal S}(\cald,J)$. If 
$\cals =(\xi,\eta,\Phi,g)\in {\mathcal S}(\cald,J)$, we have that 
$\eta_0(\xi)>0$, and $\eta={\displaystyle \frac{\eta_0}{\eta_0(\xi)}}$.
\end{proposition}

{\it Proof}. Using the canonical splitting (\ref{cs}), we write any $1$-form 
$\eta$ as $\eta=f\eta_0 +\gra$, with  $\eta_0$ and $\gra$ orthogonal to each
other. As the kernels of $\eta$ and $\eta_0$ equal ${\mathcal D}$, we 
must have $\gra=0$. Since $\eta$ is a contact form, the function $f$ is 
nowhere vanishing, and since $\Phi \mid_\cald =J=\Phi_0 \mid_\cald$, $f$ must 
be positive. The result follows. \qed

We thus see that the underlying CR structure fixes both, orientation and 
co-orientation of the contact structure. 

\begin{definition}\label{poscr} 
Let $(\cald,J)$ be a strictly pseudoconvex {\rm CR} structure of Sasaki type.
We say that a vector field $X\in \gc\gr(\cald,J)$ is {\bf positive} if 
$\eta(X)>0$ for any $\cals=(\xi,\eta,\Phi,g)\in  {\mathcal S}(\cald,J)$. 
We denote by $\gc\gr^+(\cald,J)$ the subset of all positive elements of 
$\gc\gr(\cald,J)$.
\end{definition}

We consider the mapping $\iota$ defined by projection,
\begin{equation}
\begin{array}{ccc}
{\mathcal S}(\cald,J) & \stackrel{\iota}{\rightarrow} &
 \gc\gr^{+}(\cald,J) \\ \cals & \mapsto & \xi
\end{array} \, .\label{ide}
\end{equation}
By Proposition \ref{poseta}, we see that this mapping is injective. 

We have the following.

\begin{lemma}\label{cr+lemma}
Let $(\cald,J)$ be a strictly pseudoconvex {\rm CR} structure of Sasaki type. 
Then 
\begin{enumerate}
\item[a)] $\gc\gr^+(\cald,J)$ is naturally identified with 
${\mathcal S}(\cald,J)$,
\item[b)] $\gc\gr^+(\cald,J)$ is an open convex cone in $\gc\gr(\cald,J)$,
\item[c)] The subset $\gc\gr^+(\cald,J)$ is invariant under the adjoint 
action of the Lie group $\gC\gR(\cald,J)$.
\end{enumerate}
\end{lemma} 

{\it Proof}. In order to prove (a), we show that the map $\iota$ in (\ref{ide})
is surjective. As in Proposition \ref{poseta}, we fix a Sasaki 
structure 
$\cals_0=(\xi_0,\eta_0,\Phi_0,g_0)$
in ${\mathcal S}(\cald,J)$. For $\xi\in \gc\gr^+(\cald,J)$, we define a 
$1$-form $\eta$ by 
$$\eta= \frac{\eta_0}{\eta_0(\xi)}\, .$$
Then, $\eta(\xi)=1$, and since $\xi\in \gc\gr^+(\cald,J)$, $\xi$ leaves $\cald$
invariant. This implies that $\xi \hok d\eta=\pounds_\xi\eta=0$. Thus, 
$\xi$ is the Reeb vector field of $\eta$. We then define 
$\Phi$ by $\Phi=\Phi_0-\Phi_0(\xi)\otimes \eta$, and a metric
$g$ by (\ref{me1}). The structure $\cals =(\xi,\eta,\Phi,g)$ belongs to 
${\mathcal S}(\cald,J)$, and thus, $\iota$ is surjective.

For the proof of (b), we observe that $\gc\gr^+(\cald,J)$ is open and 
convex, and that if $\xi\in \gc\gr^+(\cald,J)$, then so is $a\xi$ for any 
positive real number $a$. Indeed, all of these follow by the defining 
condition of positivity of a vector field $X$ in $\gc\gr(\cald,J)$.

For the final assertion, we observe that for groups of transformations,
the adjoint action is that induced by the differential. Thus, given 
$\phi\in \gC\gR(\cald,J)$ and $\xi \in \gc\gr^+(\cald,J)$,
we have that $\eta_0(\phi_*\xi)=(\phi^*\eta_0)(\xi)=
f_\phi\eta_0(\xi)>0$ for some positive function $f_\phi$, which
shows that $\phi_{*}\xi \in \gc\gr^+(\cald,J)$. Thus,
$\gc\gr^+(\cald,J)$ is invariant. \qed

Hereafter, we shall identify the spaces ${\mathcal S}(\cald,J)$ and 
$\gc\gr^+(\cald,J)$. We are interested in the action of the Lie group 
$\gC\gR(\cald,J)$ on ${\mathcal S}(\cald,J)=\gc\gr^+(\cald,J)$.

\begin{theorem}\label{Sascone}
Let $M$ be a closed manifold of dimension $2n+1$, and let $(\cald,J)$ be 
a {\rm CR} structure of Sasaki type on it. Then the Lie 
algebra $\gc\gr(\cald,J)$ decomposes as $\gc\gr(\cald,J)= \gt_k+ \gp$,
 where $\gt_k$ is the Lie algebra of a maximal torus $T_k$ of dimension $k$,
$1\leq k\leq n+1$, and $\gp$ is a completely reducible 
$T_k$-module. Furthermore, every $X\in \gc\gr^+(\cald,J)$ is conjugate to a 
positive element in the Lie algebra $\gt_k$. 
\end{theorem}

{\it Proof}. Let us assume first that $M$ is not the sphere with its 
standard CR structure. By Proposition \ref{autSas}, there is a Sasaki 
structure $\cals_0\in  {\mathcal S}(\cald,J)$ such that 
$\gC\gR(\cald,J)=\gA\gu\gt(\cals_0)$, which is a compact Lie group. A well 
known Lie theory result implies that every 
element in the Lie algebra $\ga\gu\gt(\cals_0)$ is conjugate under the 
adjoint action of the group $\gA\gu\gt(\cals_0)$ to one on $\gt_k$,
and by (3) of Lemma \ref{cr+lemma}, the positivity is preserved under this 
action. The possible restriction on the dimension of the maximal torus of
$\gc\gr(\cald,J)$ is well-known in Sasakian geometry.

In the case where $(\cald,J)$ is the standard CR structure on the sphere, we 
know \cite{web77} that $\gC\gR(\cald,J)={\mathbb S}{\mathbb U}(n+1,1)$, and
$\gc\gr(\cald,J)=\gs\gu(n+1,1)$, which has several maximal Abelian subalgebras.
A case by case analysis shows that the only Abelian subgroup 
where the positivity condition can be satisfied is in that of a maximal torus. 
(This can be ascertained, for instance, by looking at Theorem 6 of 
\cite{Dav05}.) \qed

We wish to study further the action of the Lie group $\gC\gR(\cald,J)$ on the 
space 
${\mathcal 
S}(\cald,J)$. The isotropy subgroup of an element 
$\cals\in {\mathcal S}(\cald, J)$ is, by 
definition, $\gA\gu\gt(\cals)$, and this contains the torus $T_k$. 
More generally, we have

\begin{lemma}\label{autS}
Let $(\cald,J)$ be a {\rm CR} structure of Sasaki type on $M$. For each 
$\cals\in {\mathcal S}(\cald, J)$, the isotropy subgroup of 
$\gC\gR(\cald,J)$ at $\cals$ is precisely $\gA\gu\gt(\cals)$.  
Furthermore,
$$\bigcap_{\cals \in {\mathcal S}(\cald, J)} \gA\gu\gt(\cals)= T_k\, .$$
In particular, $T_k$ is contained in the isotropy subgroup of every 
$\cals\in  {\mathcal S}(\cald, J)$.
\end{lemma}

{\it Proof}. It suffices to show that for the generic Reeb vector field 
$\gA\gu\gt(\cals)= T_k$. So let $\xi \in \gc\gr^+(\cald,J)$ be such that 
the leaf closure of $\calf_\xi$ is a $k$-dimensional torus $T_k$. 
Since the Reeb field is in the center of $\gA\gu\gt(\cals)$, continuity 
implies that all of $T_k$ is in the center of $\gA\gu\gt(\cals)$. But since 
the center is Abelian and $T_k$ is maximal, the result follows. \qed

We are interested in the orbit space ${\mathcal S}(\cald,J)/\gC\gR(\cald,J)$.
We have

\begin{definition}
Let $(\cald,J)$ be a {\rm CR} structure of Sasaki type on $M$.
We define the {\bf Sasaki cone} $\kappa(\cald,J)$ to be the moduli space of 
Sasakian structures compatible with $(\cald,J)$,
$$\kappa(\cald,J)={\mathcal S}(\cald,J)/\gC\gR(\cald,J) \, .$$
\end{definition}

Theorem \ref{Sascone} together with the mapping {\rm (\ref{ide})} says that 
each orbit can be 
represented by choosing a positive element in the Lie algebra $\gt_k$ of a 
maximal torus $T_k$. We denote the subset of positive elements by $\gt_k^+$,
so we have an identification $\gt_k^+=\gt_k\cap \gc\gr^+(\cald,J)\approx 
\kappa(\cald,J)$.

Now the basic Chern class of a Sasakian structure $\cals=(\xi,\eta,\Phi, g)$ 
is represented by the Ricci form $\rho^T/2\pi$ of the transverse metric 
$g_T$ (up to a factor of $2\pi$). Although the notion of basic changes with 
the Reeb vector field, the complex vector bundle $\cald$ remains fixed. Hence,
for any Sasakian structure $\cals\in {\mathcal S}(\cald,J)$, the transverse 
$2$-form $\rho^T/2\pi$ associated to $\cals$ represents the first Chern class 
$c_1(\cald)$ of the complex vector bundle $\cald$.

It is of interest to consider the case where $k=1$, that is to say, the case
where the maximal torus of $\gC\gR(\cald,J)$ is one dimensional. Since the 
Reeb vector field is central, the hypothesis that $k=1$ implies that 
$\dim \ga\gu\gt(\cals)= \dim \gc\gr(\cald,J) =1$. Hence, we have that 
${\mathcal S}(\cald,J)=\gc\gr^+(\cald,J)={\mathfrak t}_1^+=\bbr^+$, and 
${\mathcal S}(\cald,J)$ consists of the $1$-parameter family of Sasaki 
structures given by $\cals_a=(\xi_a,\eta_a,\Phi_a,g_a)$, where 
\begin{equation}
\xi_a=a^{-1}\xi\, , \quad \eta_a=a\eta\, , \quad \Phi_a=\Phi\, , \quad 
g_a=ag+(a^2-a)\eta\otimes \eta\, , \label{homo}
\end{equation}
and $a\in \bbr^+$, the $1$-parameter family of transverse homotheties.

In effect, the homotheties described above are the only deformations 
$(\xi_t,\eta_t,\Phi_t,g_t)$ of a given structure $\cals=(\xi,\eta,\Phi,g)$ in 
the Sasaki cone $\kappa(\cald,J)$ where the Reeb vector field varies in the 
form $\xi_t=f_t \xi$, $f_t$ a scalar function. For we then have that the 
family of tensors $\Phi_t$ is constant, and since $\pounds_{\xi_t} \Phi_t=0$, 
we see that 
$f_t$ must be annihilated by any section of the sub-bundle ${\mathcal D}$. 
But then (\ref{tkf}) implies that $df_t=(\xi f_t)\eta$, and we conclude that 
the function $f_t$ is constant. Thus, in describing fully the tangent 
space of ${\mathcal S}(\cald,J)$ at $\cals$, it suffices to
describe only those deformations $(\xi_t,\eta_t,\Phi_t, g_t)$ where
$\dot{\xi}= \partial_t \xi_t \mid_{t=0}$ is $g$-orthogonal to $\xi$. 
These correspond to deformations where the volume of $M$ in the 
metric $g_t$ remains constant in $t$, and are parametrized by 
elements of ${\mathfrak k}(\cald,J)$ that are $g$-orthogonal to
$\xi$.
 
The terminology we use here is chosen to emphasize the fact that the Sasaki 
cone is to a 
CR structure of Sasaki type what the K\"ahler cone is to a complex manifold 
of K\"ahler type. Indeed, for any 
point $\cals=(\xi,\eta,\Phi,g)$ in $\kappa(\cald,J)$, the complex normal bundle
$(\nu({\mathfrak F}_{\xi}),\bar{J})$ is isomorphic to $(\cald,J)$, and so is 
the underlying CR structure of any element of ${\mathcal S}(\xi,\bar{J})$.
In this sense, the complex structure $\bar{J}$ is fixed with the fixing of  
$(\cald,J)$, the Reeb vector field $\xi$ polarizes the manifold, and the 
Sasaki cone $\kappa(\cald,J)$ is the set of all possible 
polarizations.

\begin{definition}
We say that $(\xi,\eta,\Phi,g)\in \kappa(\cald,J)$ is a {\bf canonical} element
of the Sasaki cone if the space ${\mathcal S}(\xi,\bar{J})$ admits a 
canonical representative. We denote by ${\mathfrak e}(\cald,J)$ the 
set of all canonical elements of the Sasaki cone, and refer to it as the 
{\bf canonical Sasaki set} of the {\rm CR} structure $(\cald,J)$.
\end{definition}

By the identification of $\kappa(\cald,J)$ with 
$\gt_k^+$, the canonical Sasaki cone singles
out the subset of positive Reeb vector fields $\xi$ in $\gt_k^+$ for which
the functional (\ref{var}) admits a critical point.

\section{Openness of the canonical Sasaki set}
\setcounter{theorem}{0}
Given a canonical Sasakian structure $(\xi,\eta,\Phi,g)$ with underlying
CR structure $(\cald, J)$, its isometry group will contain the torus $T_k$ of
Theorem \ref{Sascone}. In fact by Lemma \ref{autS}, for a generic element 
$\xi\in 
\gt_k^+,$ $T_k$ 
will be exactly the isometry group of $g$. Moreover, Theorem \ref{ml} says 
that 
the isometry group of the transversal metric $g^T$ is a maximal compact 
subgroup $G$ of the identity component of 
the automorphism group of the transverse holomorphic structure, and the 
reductive part 
of the Lie 
algebra $\gh^T(\xi,\bar{J})/L_\xi$ consists of  
the complexification of the Lie algebra of all Killing vector fields for
$g^T$ that are Hamiltonian. If ${\mathfrak g}$ is the Lie algebra of $G$
and ${\mathfrak g}_0$ is the ideal of 
Killing fields of $g^T$ that have zeroes, then ${\mathfrak z}_0\subset
{\mathfrak g}_0$ consists of sections of $(\cald,J)$ that are 
Hamiltonian gradients, and these can be lifted \cite{BG06} to infinitesimal 
automorphisms of $(\xi,\eta,\Phi,g)$. 
Hence, these vector fields are the transversal gradients of functions in 
$M$ that are $T_k$-invariant, or to put it differently, they are 
generated by those elements of the Lie algebra ${\mathfrak t}_k$ that 
correspond to transversal holomorphic gradient sections of $(\cald,J)$. 
Thus, in searching for canonical representatives of elements of the Sasaki 
cone $\kappa(\cald,J)$, it will suffice to consider Sasakian structures that 
are invariant under $T_k$, and then seek the canonical representatives among 
them. 

We denote by ${\mathcal S} (\xi,\bar{J})^{T_k}$ the collection of all Sasakian 
structures in ${\mathcal S}(\xi,\bar{J})$ that are $T_k$-invariant,
and  by ${\mathfrak M}(\xi,\bar{J})^{T_k}$ the space of all $T_k$-invariant 
metrics
in ${\mathfrak M}(\xi,\bar{J})$. The observation made above indicates that,
in order to seek canonical representatives of ${\mathcal S}(\xi,\bar{J})$, it 
would suffices to do so among metrics in ${\mathfrak M}(\xi,\bar{J})^{T_k}$.

Given $(\xi,\eta,\Phi,g)\in {\mathcal S}(\xi,\bar{J})^{T_k}$, we let
$G_g$ stand for the Green's operator of $g$ acting on functions. We consider
a basis $\{ X_1, \ldots, X_{k-1}\}$ of ${\mathfrak z}_0\cap {\mathfrak t}_k$. 
Then the set of functions
\begin{equation}
\begin{array}{rcl}
p_0(g ) & = &  1 \\ p_j ( g ) & = &
2i G_g{\db}^{*}_{g}((JX_{j}+iX_{j})\hok d\eta )
\, , \quad j=1, \ldots , k-1\, , 
\end{array} \label{hp}
\end{equation}
spans the space of $T_k$-invariant basic real-holomorphy potentials, 
real-valued functions solutions of equation (\ref{lic}) whose $g^T$-gradients 
are holomorphic vector fields. Since the argument function on which
$G_g$ acts in order to define $p_j(g)$ is basic, 
we could have used the Green's operator of $g^T$ above instead of $G_g$
itself. 

\begin{definition}\label{df}
We define $\pi_g$ to be the $L^2$-projection onto the space of smooth 
real holomorphic potential functions in {\rm (\ref{hp})}.
\end{definition}

By Theorem \ref{th1}, 
$(\xi,\eta,\Phi,g)$ is a canonical representative of ${\mathcal S}
(\xi,\bar{J})$ if and only if $(1-\pi_g)s_g=0$, $s_g$ the 
scalar curvature of $g$.

We denote by $L^2_{B,l,T_k}$ the Hilbert space of of 
$T_k$-invariant basic real-valued functions of class $L^2_{l}$. We 
consider deformations $(\xi_\a, \eta_\a,\Phi_\a, g_\a)$ of 
$(\xi,\eta,\Phi,g)$ in the Sasaki cone ${\kappa}(\cald, J)$ where the Reeb 
vector field varies as $\xi_{\a}=\xi +\a$. We require that $\a$ be in a 
sufficiently small neighborhood of the origin in $\gc \gr (\cald,J)$ so that 
$\xi_\a$ remains positive. For $\varphi$ in a sufficiently small 
neighborhood of the origin in $L^2_{B,l+4,T_k}$, $l>n$, we then 
consider the deformations of $(\xi_\a, \eta_\a,\Phi_\a, g_\a)$ in
${\mathcal S}(\xi_{\a},\bar{J})$ to the Sasakian structure defined by
$$\begin{array}{ccl}
\eta_{\a,\varphi} & = & \eta_\a+ d^c \varphi \, , \\
\Phi_{\a,\varphi} & = & \Phi_\a - (\xi_\a \otimes (\eta_{\a,\varphi}-\eta_\a))
\circ \Phi_\a 
\, ,\\
g_{\a, \varphi} & =  & d\eta_{\a,\varphi}\circ (\BOne \otimes 
\Phi_{\a,\varphi} ) +\eta_{\a,\varphi}\otimes \eta_{\a,\varphi} \, .
\end{array}
$$ 
Here, for $(\a,\varphi)=(0,0)$, we have that $(\xi_\a,\eta_{\a,\varphi},
\Phi_{\a,\varphi},g_{\a, \varphi})=(\xi,\eta,\Phi,g)$.  
The restriction on $l$ ensures that the curvature 
tensors of $g_{\a,\varphi}$ are all well-defined because, under such a
constraint, $L^2_{B,l,T_k}$ is a Banach algebra. 

We let ${\mathcal U}\subset \gc \gr (\cald,J)\times L_{B,l+4,T_k}^{2}(M)$ be 
the open neighborhood of $(0,0)$ where the two-parameter family of deformations
$g_{\a,\v}$ of $g$ is well-defined, and consider the scalar curvature map
\begin{equation}
\begin{array}{cccc}
\gc \gr (\cald,J) \times L_{B,l+4,T_k}^{2}(M) \supset &{\mathcal U} &
\stackrel{{S}}{\longrightarrow }& L_{B,l,T_k}^{2}(M) \\
&(\a , \varphi )
&\mapsto & s_{g_{\a,\varphi}}\, ,
\end{array}\label{sca1}
\end{equation}
where $s_{g_{\a,\varphi}}$ is the scalar curvature of the metric 
$g_{\a,\varphi}$. 

\begin{proposition}
For $l>n$, the map {\rm (\ref{sca1})} is well-defined and $C^{1}$, with
Fr\'echet derivative at the origin given by
\begin{equation}
DS_{(0,0)}=\left[ -n \Delta_B \left( \eta( \, 
\cdot \, )\right)  \hspace{.2in}
-\frac{1}{2}(\Delta_{B}^2 + 2r^T \cdot \nabla_{T} \nabla _{T})
\right] \, , \label{va1}
\end{equation}
where the quantities in the right are associated to the transversal metric
$g^T$ defined by $g$, and $r^{T} \cdot \nabla_{T} \nabla _{T}$ denotes the 
full contraction of the Ricci tensor and two covariant derivatives of $g^T$.
\end{proposition}

{\it Proof}. Notice that when deforming $(\xi,\eta,\Phi,g)$ to
$(\tilde{\xi},\tilde{\eta},\tilde{\Phi},\tilde{g})$ while preserving
the underlying CR structure $(\cald,J)$, the transversal K\"ahler form
$d\eta$ changes by the conformal factor $f(\tilde{\xi})=1/\eta(\tilde{\xi})$.
The first component of the Fr\'echet derivative above follows via a simple
calculation, after observing that
the Ricci tensor of the transversal K\"ahler metric is computed in 
a holomorphic frame by $-i\ddb \log \det{(g_{i\bar{k}}^T)}$.
The second component of the Fr\'echet derivative follows by
Proposition \ref{pr2}.  \qed
\medskip

For  any integer $l$, we let $I_{l}\subset L^2_{B,l ,T_k}$ denote the 
orthogonal complement of the kernel of the operator 
$L^B_g=(\bar{\d}{\d}_{g}^{\#})^{\ast}\bar{\d}{\d}_{g}$ in (\ref{li}), 
and set ${\mathcal V}= {\mathcal U}
\cap (\gc \gr (\cald, J)\times I_{k+4})$, where ${\mathcal U}$ is the 
neighborhood of $(0,0)$ in $\gc \gr (\cald, J)\times L^2_{B,k+4,T_k}$ in
(\ref{sca1}), shrunk if necessary so that 
$${\rm ker}(1-\pi_g)(1-\pi_{g_{\a,\varphi}})={\rm ker}
(1-\pi_{g_{\a,\varphi}})$$
whenever  we have a Sasaki metric $g_{\a,\varphi}$ 
of the type indicated above, parametrized by some $(\a ,\varphi )\in 
{\mathcal U}$.
Here, $\pi_g$ is the projection onto the finite dimensional space of 
functions (\ref{hp}) introduced in Definition \ref{df}. It is clear that
the range of this projection changes smoothly with the metric.

Since a Sasaki metric $g$ is canonical if, and only if, its scalar curvature 
is annihilated by the projection operator $1-{\pi}_{g}$,
we introduce the map
\begin{equation}
\begin{array}{rcl}
\gc \gr (\cald, J)\times I_{l+4} \supset \, {\mathcal V} & 
\stackrel{{\mathcal S}}{\longrightarrow} &  \gc \gr (\cald, J)\times I_{k} \\
{\mathcal S}(\a,\varphi ) & := & (\a ,(1-\pi_{g})(1-{\pi}_{g_{\a,\varphi}})
\, S(\a , \varphi  )\, ) 
\end{array}\, , \label{om}
\end{equation}
where $S(\a , \varphi  )$ is the map in (\ref{sca1}).

We have the following.

\begin{lemma}
Suppose that $g$ is a canonical Sasaki metric representing the polarization of
$M$ given by $(\xi,\bar{J})$, and that $g_{t}=g_{t,\a,\varphi}$ is a curve of 
Sasakian metrics of the type above that starts at $g$ when $t=0$, and is
parametrized by $(\a,\varphi)\in \gc \gr (\cald,J)\times L_{B,l+4,T_k}^2$.
Then
$$(1-\pi_{g})\! \left.(\frac{d}{dt}\pi _{g_{t}})\right|_{t=0}\! \! \! s_g
=(1-\pi_{g})
[ -2iG_g{\db}^{*}\left( \left(  \eta(\a)\right)\db s_g 
\right) + 
({\d}s_g \hok {\d}^{\#} \varphi )]\, ,$$
where $G_g$ is the Green's operator of $g$.
\end{lemma}

{\it Proof}. The result is clear if the basic scalar curvature function $s_g$
is constant. For the general case, we refer the reader to the 
original argument in \cite{cs} for the K\"ahler case (see also the 
significantly improved  understanding, and related discussions, given in 
\cite{si2}). Its required extension follows by the observation that the
transversal metric of $g$ is deformed by the conformal factor 
$1/\eta({\xi}_{t,\a,\v})$, as noted before.
\qed
\medskip

\begin{proposition}\label{ipro}
For $l>n$, the map {\rm (\ref{om})} is $C^{1}$ with 
Fr\'echet derivative at the origin given by
\begin{equation}
D{\mathcal S}_{(0,0)}\! =\! \left(\! \!  \begin{array}{cc} 1 & 0 \\ 0 &
1\! -\! \pi_{g} \end{array} \! \! \right) \left( \! \! \begin{array}{cc}
1 & 0 \\ -n\Delta_B\left( \eta( \, \cdot \, )\right)
\!+\! 2 iG_g{\db}^{*} \! \! \left( \left( \! \eta( \, \cdot \, )
\right)
\db s_g \! \right) \! &  -2L_g^B
\end{array} \! \! \right) \, , \label{si3}
\end{equation}
where $L_g^B= (\bar{\d}{\d}^{\#})^{\ast}\bar{\d}{\d}^{\#}$.
\end{proposition}

\begin{proposition}
Let $M$ be a closed manifold, and $(\cald,J)$ be a {\rm CR} structure of
Sasaki type on it. Then the
map ${\mathcal S}$ defined in {\rm (\ref{om})} becomes a diffeomorphism when
restricted to a sufficiently small neighborhood of the origin.
\end{proposition}

{\it Proof}. We apply the inverse function theorem for Banach spaces. Hence,
we just need to prove that $D{\mathcal S}_{(0,0)}$ has trivial kernel and
cokernel.

Suppose that $(\a,\v)$ is in the kernel of $D{\mathcal S}_{(0,0)}$. By
(\ref{si3}), we see that $\a=0$, and that $$(1-\pi_{g}) L_g^B \v =0\, .$$
It follows that $L_g^B\v=(\bar{\d}{\d}^{\#})^{\ast}\bar{\d}{\d}^{\#}\v$ is a 
holomorphy potential, and consequently, it can be written as
$$(\bar{\d}{\d}^{\#})^{\ast}\bar{\d}{\d}^{\#}\v =\sum_j c_j f_{g}^j \, ,$$
in terms of an orthonormal basis of the space spanned by the functions
in (\ref{hp}). If we take the inner product of this expression with $f_g^j$,
and dualize the symmetric map $(\bar{\d}{\d}^{\#})^{\ast}\bar{\d}{\d}^{\#}$,
we see that $c_j=0$. Thus, $(\bar{\d}{\d}^{\#})^{\ast}\bar{\d}{\d}^{\#}\v =0$.
But $\v \in I_{l+4}$, space orthogonal to the kernel of
 $(\bar{\d}{\d}^{\#})^{\ast}\bar{\d}{\d}^{\#}$. So $\v$ must be zero, and the
kernel of $D{\mathcal S}_{(0,0)}$ consists of the point $(0,0)$.

Suppose now that $(\beta, \psi)$ is orthogonal to every element in the image
of $D{\mathcal S}_{(0,0)}$. Then, it must be orthogonal to the image of
$(0,\v)$ for any $\v\in I_{l+4}$, and therefore,
$$\< (1-\pi_g)(\bar{\d}{\d}^{\#})^{\ast}\bar{\d}{\d}^{\#}\v,\psi\>=
\< \v, (\bar{\d}{\d}^{\#})^{\ast}\bar{\d}{\d}^{\#}(1-\pi_g)\psi\> =0$$
for all such $\v$. It follows that the component of
$(\bar{\d}{\d}^{\#})^{\ast}\bar{\d}{\d}^{\#}(1-\pi_g)\psi$ perpendicular
to the kernel of $(\bar{\d}{\d}^{\#})^{\ast}\bar{\d}{\d}^{\#}$ is zero, and
thus, $(\bar{\d}{\d}^{\#})^{\ast}\bar{\d}{\d}^{\#}(1-\pi_g)\psi=\sum c_j
f_g^j$. The same argument used above implies that $c_j=0$, and so,
$(1-\pi_g)\psi$ is the kernel of $(\bar{\d}{\d}^{\#})^{\ast}\bar{
\d}{\d}^{\#}$. But the image of $1-\pi_g$ is orthogonal to this kernel.
Hence, $\psi=0$.  
Using this fact, we may now conclude that $\beta$ must be such that
$\< \beta, \a\> =0$ for all $\a \in \gc \gr (\cald,J)$, and so
$\b=0$. The cokernel of $D{\mathcal S}_{(0,0)}$ is
trivial. \qed

We now prove the following result, analogous to the openness of the extremal
cone in K\"ahler geometry \cite{cs}. 

\begin{theorem}\label{ss}
Let $(\cald,J)$ be a {\rm CR} structure of Sasaki type on $M$. Then the 
canonical Sasaki set ${\mathfrak e}(\cald,J)$ is an open subset of 
the Sasaki cone $\kappa(\cald,J)$.
\end{theorem}

{\it Proof}. Let ${\mathcal V}_0\subset {\mathcal V}$ be a neighborhood of 
$(0,0) \in \gc \gr (\cald,J)\times I_{k+4}$ such
that ${\mathcal S}\mid_{{\mathcal V}_0}$ is a diffeomorphism from 
${\mathcal V}_0$ onto an open neighborhood of the origin in 
$\gc \gr (\cald, J)\times I_{k}$. For any point $\a$ in
${\mathcal S}({\mathcal V}_0) \cap ( \gc \gr (\cald, J)\times \{0\} )$, we
define $\varphi(\a)$ to be the projection onto $I_{k+4}$ of
$({\mathcal S}\mid_{{\mathcal V}_0})^{-1}(\a) $. Then we have
$$(\a,0)={\mathcal S}(\a,\varphi(\a))=(\a, (1-\pi_{g})(1-\pi_{g_{\a,\v(\a)}})
s_{g_{\a,\varphi(\a)}}) \, ,$$
where $g_{\a,\v(\a)}$ is the deformation of the metric $g$ associated to the
parameters $(\a,\v(\a))$. Since the kernel of
$(1-\pi_{g})(1-\pi_{g_{\a,\v(\a)}})$ equals the kernel of
$1-\pi_{g_{\a,\v(\a)}}$, it follows that
$(1-\pi_{g_{\a,\v(\a)}})s_{g_{\a,\varphi(\a)}}=0$. We then have that the
scalar curvature of the Sasaki metric $g_{\a,\v(\a)}$ is a holomorphy 
potential, and so this metric is a canonical representative of
${\mathcal S}(\xi+\a, \bar{J})$. This completes the proof. \qed

\begin{example}
Let us take coordinates $z=(z_0, \ldots , z_{n})$ in  ${\mathbb C}^{n+1}$,
and consider the unit sphere
$${\mathbb S}^{2n+1}=\{ z \in {\mathbb C}^{n+1}: \; |z|=1 \} \, .$$
If $z_k=x_k+iy_k$ is the decomposition of $z_k$ into real and imaginary parts, 
then
the vector fields $H_k = (y_k  \partial_{x_k} -x_k \partial_{y_k})$,
$k=0,\ldots, n$, form a basis for the Lie algebra ${\mathfrak t}_{n+1}=
{\mathbb R}^{n+1}$ of a maximal torus in the automorphism group
${\mathbb U}(n+1)$ of the standard Sasakian structure on
${\mathbb S}^{2n+1}$. The latter is given by the contact form
$\eta= \sum_{k=0}^n (y_k dx_k -x_k dy_k)$, Reeb vector field
$\xi= \sum_{k=0}^n H_k $, and $(1,1)$-tensor $\Phi$ defined by the
restriction $J$ of the complex structure on ${\mathbb C}^{n+1}$ to 
${\mathcal D}={\rm ker}\, \eta$, and the fact that $\Phi (\xi)=0$. The 
ensuing compatible metric (\ref{me1}) defined by $(\xi,\eta,\Phi)$ is the
standard metric $g$ on ${\mathbb S}^{2n+1}$, which is Einstein with
$r_g=2n g$. Thus, the metric $g$ yields a canonical representative of 
${\mathcal S}(\xi, J)$. We have that the
automorphism group of $(\xi,\eta,\Phi,g)$ is $\bbu(n+1)$. Its maximal 
torus $\gT_{n+1}$ has Lie algebra $\gt_{n+1}$ with basis 
$\{ H_0, \ldots, H_n\}$.

We now fix this CR structure $({\mathcal D},J)$ on ${\mathbb S}^{2n+1}$, 
and consider the set ${\mathcal S}({\mathcal D},J)$ of all Sasakian structures
associated with it. A vector field $X$ is positive if,
and only if, $\eta(X)>0$, and such a vector field is conjugate to a positive
vector in the Lie algebra $\gt_{n+1}$. We use the basis 
$\{ H_0,\ldots, H_n\}$ to identify this Lie algebra with $\bbr^{n+1}$, so the 
point $w=(w_0,\ldots, w_n)$ in $\bbr^{n+1}$ yields the vector 
$\xi_w=\sum w_k H_k$. Then we have that
$$\eta(\xi_w)= \sum_{i=0}^n w_i\bigl(x_i^2+y_i^2\bigr)=
\sum_{i=0}^nw_i|z_i|^2\, ,
$$
and so, the set of positive elements of $\gt_{n+1}$ is just $\bbr_{+}^{n+1}$.
By Theorem \ref{Sascone}, the Sasaki cone 
$\kappa(\cald,J)$ is equal to $\bbr_{+}^{n+1}$. If 
$w\in \bbr_{+}^{n+1}$, this vector gives rise
to the Sasakian structure $(\xi_w,\eta/\eta(\xi_w),\Phi_w,g_w)$,
where $\Phi_w$ is defined by the conditions $\Phi_w\mid_{\cald}=J$ and 
$\Phi_w(\xi_w)=0$, respectively, and $g_w$ is 
determined by the expression (\ref{me1}) in terms of $\xi_w$, 
$\eta/\eta(\xi_w)$ and $\Phi_w$.

For any $w \in {\mathbb Z}_{+}^{n+1}$, the 
Sasakian structure $(\xi_w,\eta/\eta(\xi_w),\Phi_w,g_w)$ is quasi-regular, and 
its transversal is a manifold with orbifold singularities, the weighted 
projective space ${\mathbb C}{\mathbb P}^n_w$. The space of metrics 
${\mathfrak M}(\xi_w,
\bar{J})$ associated with the polarized Sasakian manifold 
$({\mathbb S}^{2n+1},\xi_w, \bar{J})$ has a representative $g_w$ whose
transverse K\"ahler metric $g_w^T$ is Bochner flat \cite{br} on
${\mathbb C}{\mathbb P}^n_w$, and thus, extremal. Computing in an affine 
orbifold chart, it can be determined \cite{dapa} that the scalar curvature
of $g_w^T$ is given by
$$s_{g_w^T}=4(n+1)\frac{\sum_{j=0}^n w_j(2(\sum_{k=0}^n w_k) -(n+2)w_j)
|z_j|^2 }{\sum_{j=0}^n w_j |z_j|^2 }\, ,$$
at $z\in {\mathbb S}^{2n+1}$. Since the volume $\mu_{g_w^T}({\mathbb C}
{\mathbb P}^n_w) =\pi^n/(n! \prod_{j=0}^n w_j )$ \cite{dapa}, 
the volume of ${\mathbb S}^{2n+1}$ in the Sasakian metric $g_w$ is just
$$\mu_{g_w}({\mathbb S}^{2n+1})= 2 \frac{\pi^{n+1}}{n!}\frac{1}
{\prod_{j=0}^n w_j}\, .$$
Similarly, since $s_{g_w}=s_{g_w^T}-2n$, and since the mean transverse scalar 
curvature of $g_w^T$ is $4n\sum _{j=0}^n w_j$ \cite{dapa}, we have that the 
projection $s^0_{g_w}$ of $s_{g_w}$ onto the constants is given by
$$s^0_{g_w}=2n(2 \sum_{j=0}^n w_j -1)\, .$$
Thus,
$$s_{g_w}-s^0_{g_w}=4(n+2)\frac{\sum_{j=0}^n w_j((\sum_{k=0}^n w_k) -(n+1)w_j)
|z_j|^2 }{\sum_{j=0}^n w_j |z_j|^2 }\, ,$$
Notice that if the weight vector $w$ is of the
form $w=l(1, \ldots, 1)$, then
$s_{g_w^T}=4(n+1)n l$, and this yields the 
scalar curvature of the Fubini-Study metric when $l=1$, as it should. 

For the Sasaki-Futaki character of the polarization $(\xi_w,\bar{J})$, it 
suffices to determine its value on vector fields $X$ that commute
with $\xi_w$, and that are of the form $X=\d^{\#}f$ for $f$ a basic real 
holomorphy potential. In that case, we have
$${\mathfrak F}_{(\xi_w,\bar{J})}(X)=-\int_{{\mathbb S}^{2n+1}}f (s_{g_w}-
s_{g_w^0}) d\mu_{g_w} \, .$$
For convenience, let us set $A_j=\sum_{k=0}^n w_k -(n+1)w_j$. Working in an 
affine orbifold chart for ${\mathbb C}{\mathbb P}_w$, we 
then see that if $f=\sum_{i=0}^n b_i |z_i|^2$ we have that
$$\begin{array}{rcl}
{\mathfrak F}_{(\xi_w,\bar{J})}(X) & = &{\displaystyle -8(n+2)\pi ^{n+1} 
\int_{{\mathbb R}^{n}_{+}}
\frac{(b_0+\sum_{j=1}^n b_jx_j)(w_0A_0+\sum_{j=1}^n w_jA_jx_j )}
{(w_0+\sum_{j=1}^n w_j x_j)^{n+3} }dx_1\ldots dx_n } \vspace{1mm} \\
& = & {\displaystyle -16\frac{\pi ^{n+1}}{(n+1)!}\left(
\sum_{i=0}^n \frac{b_i}{w_i}A_i +\frac{1}{2}\sum_{i\neq j}\frac{b_i}{w_i}A_j
\right) \frac{1}{\prod_{j=0}^n w_j}}\, .
\end{array}$$

\begin{theorem}
Let $(\cald,J)$ be the standard {\rm CR} structure on the unit sphere
${\mathbb S}^{2n+1}$,
and $\kappa(\cald,J)$ and ${\mathfrak e}(\cald,J)$ be the associated Sasaki 
cone and canonical Sasaki set, respectively. Then 
${\mathfrak e}(\cald,J)= \kappa(\cald,J)$, and the only canonical points in 
${\mathfrak e}(\cald,J)$ that yield metrics of constant scalar curvature are 
those representing transverse homotheties of the standard Riemannian Hopf 
fibration. Of these, the metric of constant sectional curvature one
is the only Sasaki-Einstein metric whose underlying {\rm CR} structure 
is $(\cald, J)$.
\end{theorem}

{\it Proof}. We have identified above $\kappa(\cald,J)$ with points $w$ in 
${\mathbb R}_{+}^{n+1}$. If $w\in {\mathbb Z}_{+}^{n+1}$, the polarized 
Sasakian manifold $({\mathbb S}^{2n+1},\xi_w, \bar{J})$ admits a representative
$g_w$ with transverse Bochner flat metric on the transverse space
${\mathbb C}{\mathbb P}^n_w$. Thus, these weights $w$ belong to
${\mathfrak e}(\cald,J)$. Using the homotheties 
(\ref{homo}), we may obtain canonical representatives of the 
polarized Sasakian manifold $({\mathbb S}^{2n+1},\xi_w, \bar{J})$ for
any weight $w\in {\mathbb Q}_{+}^{n+1}$. Applying Theorem
\ref{ss}, we obtain the same result for arbitrary weights $w$ in
${\mathbb R}_{+}^{n+1}$. Thus, ${\mathfrak e}(\cald,J)= \kappa(\cald,J)$.

The expressions computed above for the scalar curvature, volume, and 
Sasaki-Futaki character for $w=(w_0,\ldots, w_n)\in {\mathbb Z}_{+}^{n+1}$ are
rational functions of the weights $w_j$, so they also define the 
scalar curvature, volume, and Sasaki-Futaki character 
when $w$ is an arbitrary vector in ${\mathbb R}_+^{n+1}$, regardless of the
fact that there may not be a transversal manifold to speak of in this general
situation.

The assertion about the scalar curvatures of the canonical representatives
$g_w$ follows by the expression for ${\mathfrak F}_{(\xi_w, \bar{J})}$
given above, and Proposition \ref{FSscalar}. Indeed, the character  
${\mathfrak F}_{(\xi_w, \bar{J})}$ is identically zero if, and only if, 
$A_j=A_k$ for all pairs of indices $j,k$. This only happens if the vector
of weights $w$ is of the form $w=l(1,\ldots, 1)$.
The only one of these metrics that is Sasaki-Einstein is the standard
metric. This follows by a simple analysis of the change of the Ricci curvature
under homotheties of the metric.
\qed
\end{example}

In this example, the first Chern class of the sub-bundle $\cald$ is trivial.
Thus, for any $w$ in the Sasaki cone, the basic first Chern class of the 
resulting foliated manifold is proportional to the basic class defined by 
the transversal K\"ahler form $d\eta_w$. However, only for the weight
$w=(1, \ldots, 1)$ there exists a Sasaki-Einstein representative.

In general, given a CR structure $(\cald, J)$ of Sasaki type on a closed 
manifold $M$, we do not expect the equality ${\mathfrak e}(\cald,J)=
\kappa(\cald,J)$ to hold, though this is likely to be so in the toric case.


\begin{thebibliography}{99}
\bibitem{BG00}
C.P. Boyer \& K. Galicki, {\it A note on toric contact geometry}, J. Geom.
Phys., 35 (2000), pp. 288-298. 
\bibitem{BG01a}
\bysame, {\it Einstein manifolds and contact geometry}, Proc. Amer. Math. 
129 (2001), pp. 2419-2430.
\bibitem{BG06}
\bysame, {\it Sasakian Geometry}, Oxford Mathematical
Monographs, Oxford University Press, to appear, Oxford, 2006.
\bibitem{BoGa05a}
\bysame, \emph{Sasakian geometry, hypersurface singularities, and Einstein
metrics}, Rend. Circ. Mat. Palermo (2) Suppl. 75 (2005), pp. 57-87.
\bibitem{BGK05}
C.P. Boyer, K. Galicki \& J. Koll\'ar, {\it Einstein metrics on spheres},
Ann. of Math., 162 (2005), pp. 557-580.
\bibitem{BGM06}
C.P. Boyer, K. Galicki \& P. Matzeu, {\it On Eta-Einstein Sasakian
  Geometry}, Comm. Math. Phys., 262 (2006), pp. 177-208.
\bibitem{brsl}
V. Brinzanescu \& R. Slobodeanu, {\it Holomorphicity and Walczak formula
on Sasakian manifolds}, preprint 2005. 
\bibitem{br}
R. L. Bryant, Bochner-K\"ahler metrics, J. Amer. Math. Soc., 14 (2001), pp.
623-715. 
\bibitem{cal1} 
E. Calabi, {\it Extremal K{\"a}hler metrics}, in {\it
Seminar of Differerential Geometry}, ed. S. T. Yau, Annals of
Math. Studies, 102, Princeton University Press (1982), pp. 259-290.
\bibitem{cal2}
\bysame, {\it Extremal K\"ahler metrics II}, in
Differential geometry and complex analysis (I. Chavel \&
H.M. Farkas eds.), Springer-Verlag, 1985, pp. 95-114.
\bibitem{cal0}
\bysame, {\it The Space of K\"ahler Metrics}, Proc. Int. Cong. Math.,
Amsterdam, (1954), pp. 206-207.  
\bibitem{ChTi94}
J. Cheeger \& G. Tian, {\it On the cone structure at infinity of Ricci flat
manifolds with Euclidean volume growth and quadratic curvature decay},
Invent. Math., 118 (1994), pp. 493-571. 
\bibitem{ChMo74}
S.S. Chern \& J.K. Moser, {\it Real hypersurfaces in complex manifolds},
  Acta Math. 133 (1974), pp. 219-271. 
\bibitem{CLPP05}
M. Cvetic, H. L{\"u}, Don N. Page \& C.N. Pope, {\it New Einstein-Sasaki 
spaces in five and higher dimensions}, Phys. Rev. Lett., 95 (2005), 
p. 4, \MR{2167018}.
\bibitem{Dav05}
L. David, {\it The Bochner-flat cone of a CR manifold}, math.DG/0512604 (2005).
\bibitem{dapa}
L. David \& P. Gauduchon, {\it The Bochner-flat geometry of weighted projective
spaces}, C.R.M. Proceedings and Lecture Notes, 40, 2005.
\bibitem{fu}
A. Futaki, {\it An obstruction to the existence of Einstein K\"ahler
metrics}, Invent. Math., 73 (1983), pp. 437-443.  
\bibitem{GMSW04a}
J.P. Gauntlett, D. Martelli, J. Sparks \& W. Waldram, {\it Sasaki-Einstein 
metrics on $S^2\times S^3$}, Adv. Theor. Math. Phys., 8 (2004), pp. 711-734.
\bibitem{GMSW04b}
\bysame, {\it A new infinite class of Sasaki-Einstein manifolds}, 2004,
 Adv. Theor. Math. Phys., 8 (2004), pp. 987-1000.
\bibitem{Kol05a}
J. Koll\'ar, {\it Circle actions on simply connected $5$-manifolds},
arXiv:math.GT/ 0505343 (2005).
\bibitem{Kol05b}
\bysame, {\it Einstein metrics on five-dimensional Seifert bundles}, J. Geom. 
Anal., 15 (2005), pp. 445-476. 
\bibitem{cs} C. LeBrun \& S.R. Simanca, {\it On the K\"ahler Classes of 
Extremal Metrics}, Geometry and Global Analysis, (First MSJ Intern. Res. 
Inst. Sendai, Japan) eds. Kotake, Nishikawa \& Schoen, 1993.
\bibitem{Lee96}
J.M. Lee, {\it CR manifolds with noncompact connected automorphism
  groups}, J. Geom. Anal., 6 (1996), pp. 79-90.
\bibitem{Lic57}
A. Lichnerowicz, {\it Sur les transformations analytiques des
vari\'et\'es k\"ahl\'eriennes compactes}, C.R. Acad. Sci. Paris, 244
(1957), pp. 3011-3013. 
\bibitem{MaSp05b}
D. Martelli \& J. Sparks, {\it Toric Sasaki-Einstein metrics on $S^2\times 
S^3$}, Phys. Lett. B, 621 (2005), pp. 208-212, \MR{2152673}.
\bibitem{MaSp06}
\bysame, {\it Toric geometry, Sasaki-{E}instein manifolds and a new infinite 
class of {A}d{S}/{C}{F}{T} duals}, Comm. Math. Phys., 262 (2006), pp. 51-89.
\bibitem{MaSpYau05}
D. Martelli, J. Sparks \& S.-T. Yau, {\it The geometric dual of
a-maximisation for toric Sasaki-Einstein manifolds},
arXiv:hep-th/0503183, preprint 2005.
\bibitem{MaSpYau06}
\bysame, {\it Sasaki-Einstein manifolds and volume minimisation},
  arXiv:hep-th/0603021, preprint 2006.
\bibitem{Mat57}
Y. Matsushima, {\it Sur la structure du groupe d'hom\'eomorphismes
analytiques d'une certaine vari\'et\'e k\"ahl\'erienne}, Nagoya Math. J.,
11 (1957), pp. 145-150. 
\bibitem{Mol88}
P. Molino, {\it Riemannian foliations}, Progress in Mathematics 73,
Birkh\"auser Boston Inc. 1988. 
\bibitem{nito}
S. Nishikawa \& P. Tondeur, {\it Transversal infinitesimal automorphisms for
harmonic K\"ahler foliations}, T\^{o}huku Math. J., 40 (1988), pp. 599-611.
\bibitem{on}
B. O'Neill, {\it The fundamental equations of a submersion}, Mich. Math. J.,
13 (1966), pp. 459-469.
\bibitem{Sch95}
R. Schoen, {\it On the conformal and CR automorphism groups}, Geom. Funct.
Anal., 5 (1995), pp. 464-481.
\bibitem{si}
S.R. Simanca, {\it Canonical metrics on compact almost complex manifolds},
Publica\c{c}\~{o}es Matem\'aticas do IMPA, IMPA, Rio de Janeiro, 2004. 
97 pp.
\bibitem{si2}
\bysame, {\it Heat Flows for Extremal K\"ahler Metrics}, Ann. Scuola Norm.
Sup. Pisa CL. Sci., 4 (2005), pp. 187-217.
\bibitem{web77}
S.M. Webster, {\it On the transformation group of a real hypersurface},
Trans. Amer. Math. Soc., 231 (1977), pp. 179-190. 
\end{thebibliography}
\end{document}